\documentclass[VANCOUVER,STIX1COL]{WileyNJD-v2}

\usepackage{moreverb}

\usepackage{xfrac} %
\usepackage{siunitx}

\usepackage[capitalise,noabbrev]{cleveref}

\usepackage{tikz}
\usepackage{multirow}

\crefformat{equation}{(#2#1#3)}

\renewcommand{\vec}[1]{\boldsymbol{#1}}

\articletype{Original Research}%

\received{<day> <Month>, <year>}
\revised{<day> <Month>, <year>}
\accepted{<day> <Month>, <year>}

\newcommand{\ssff}{FF}
\newcommand{\ssfc}{FC}
\newcommand{\sscf}{CF}
\newcommand{\sscc}{CC}
\newcommand{\ssf}{F}
\newcommand{\ssc}{C}

\DeclareMathOperator{\argmin}{argmin}

\graphicspath{{./figures/}}

\usepackage{collcell}
\usepackage{xcolor}%
\usepackage{etoolbox}
\usepackage{pgf}
\usepackage{ifthen}
\newtoggle{inTableHeader}%
\toggletrue{inTableHeader}%
\newcommand*{\StartTableHeader}{\global\toggletrue{inTableHeader}}%
\newcommand*{\EndTableHeader}{\global\togglefalse{inTableHeader}}%

\let\OldTabular\tabular%
\let\OldEndTabular\endtabular%
\renewenvironment{tabular}{\StartTableHeader\OldTabular}{\OldEndTabular\StartTableHeader}%

\definecolor{mygreen}{RGB}{6, 156, 86}
\definecolor{myorange}{RGB}{255, 104, 30}
\definecolor{myred}{RGB}{211, 33, 44}

\newcommand*{\myColorCell}[1]{%
   \iftoggle{inTableHeader}{#1}{%
     \pgfmathsetmacro{\val}{#1}
     \pgfmathsetmacro{\thiscolor}{\val<0.4 ? "mygreen" : (\val>0.8 ? "myred" : "black")}  %
     \textcolor{\thiscolor}{#1}
   }}%
\newcolumntype{N}{>{\collectcell\myColorCell}c<{\endcollectcell}}%

\newcommand*{\myColorCellCopTwo}[1]{%
   \iftoggle{inTableHeader}{#1}{%
     \pgfmathsetmacro{\val}{#1}
     \pgfmathsetmacro{\thiscolor}{\val<1.5 ? "mygreen" : (\val<2.0 ? "black" : (\val<2.5 ? "myorange" : "myred"))}  %
     \textcolor{\thiscolor}{#1}
   }}%
\newcolumntype{\CopTwo}{>{\collectcell\myColorCellCopTwo}c<{\endcollectcell}}%

\newcommand*{\myColorCellCopThree}[1]{%
   \iftoggle{inTableHeader}{#1}{%
     \pgfmathsetmacro{\val}{#1}
     \pgfmathsetmacro{\thiscolor}{\val<2.0 ? "mygreen" : (\val<2.5 ? "black" : (\val<3.0 ? "myorange" : "myred"))}  %
     \textcolor{\thiscolor}{#1}
   }}%
\newcolumntype{\CopThree}{>{\collectcell\myColorCellCopThree}c<{\endcollectcell}}%

\newcommand*{\myColorCellCgrid}[1]{%
   \iftoggle{inTableHeader}{#1}{%
     \pgfmathsetmacro{\val}{#1}
     \pgfmathsetmacro{\thiscolor}{\val<1.5 ? "black" : (\val>1.5 ? "black" : "black")}  %
     \textcolor{\thiscolor}{#1}
   }}%
\newcolumntype{\Cgrid}{>{\collectcell\myColorCellCgrid}c<{\endcollectcell}}%

\begin{document}

\newcommand{\thistitle}{Generalizing Reduction-Based Algebraic Multigrid}
\title{\thistitle}

\author[1]{Tareq Zaman*}
\author[2]{Nicolas Nytko}
\author[3]{Ali Taghibakhshi}
\author[1]{Scott MacLachlan}
\author[2]{Luke Olson}
\author[3]{Matthew West}

\address[1]{\orgdiv{Interdisciplinary Program in Scientific Computing}, \orgname{Memorial University of Newfoundland}, \orgaddress{\state{NL}, \country{Canada}}}
\address[2]{\orgdiv{Department of Computer Science}, \orgname{University of Illinois at Urbana-Champaign}, \orgaddress{\state{Illinois}, \country{USA}}}
\address[3]{\orgdiv{Department of Mechanical Science and Engineering}, \orgname{University of Illinois at Urbana-Champaign}, \orgaddress{\state{Illinois}, \country{USA}}}

\authormark{Zaman, Nytko, Taghibakhshi, MacLachlan, Olson, and West}

\corres{*T. U. Zaman, Interdisciplinary Program in Scientific Computing, Memorial University of Newfoundland, Newfoundland and Labrador, Canada. \email{tzaman@mun.ca}}

\abstract[Abstract]{
Algebraic Multigrid (AMG) methods are often robust and effective solvers
for solving the large and sparse linear systems that arise from discretized
PDEs and other problems, relying on heuristic graph algorithms to
achieve their performance.  Reduction-based AMG (AMGr) algorithms
attempt to formalize these heuristics by providing two-level convergence bounds
that depend concretely on properties of the partitioning of the given
matrix into its fine- and coarse-grid degrees of freedom.  MacLachlan
and Saad (SISC 2007) proved that the AMGr method yields provably robust
two-level convergence for symmetric and positive-definite matrices
that are diagonally dominant, with a convergence factor bounded as a
function of a coarsening parameter.  However, when applying AMGr
algorithms to matrices that are not diagonally dominant, not only do
the convergence factor bounds not hold, but measured performance is
notably degraded.  Here, we present modifications to the
classical AMGr algorithm that improve its performance
on matrices that are not diagonally dominant, making use of strength
of connection, sparse approximate inverse (SPAI) techniques, and
interpolation truncation and rescaling, to improve robustness while
maintaining control of the algorithmic costs.
We present numerical
results demonstrating the robustness of this approach for both
classical isotropic diffusion problems and for
non-diagonally dominant systems coming from anisotropic diffusion.
}


\maketitle

\section{Introduction}\label{sec:intro-gamgr}

Partial differential equations (PDEs) arise naturally as mathematical models of physical systems in many fields of science and engineering.  As analytical techniques for their solution are limited to simple equations and geometries, numerical approximation of solutions via discretization techniques is ubiquitous.  Standard discretizations necessitate the solution of large linear and nonlinear systems of equations which, in turn, requires fast and efficient algorithms.  Among the techniques most commonly used for discretized elliptic equations are multigrid (MG) methods, known for their efficient and robust solution of a wide range of problems.  Geometric multigrid (GMG) methods are often most efficient, but require detailed knowledge of the problem to be solved, its discretization, and regular structure of the underlying mesh hierarchy.  In contrast, algebraic multigrid (AMG) methods can be effectively applied to problems on unstructured grids, or with highly variable (or discontinuous) coefficients.  While the idea of AMG was first proposed over 40 years ago~\cite{ABrandt_SFMcCormick_JWRuge_1982a, KStuben_1983a, ABrandt_SFMcCormick_JWRuge_1984a, JWRuge_KStuben_1987a}, understanding and improving the convergence of AMG remains an active area of research.

Classical (Ruge-St\"{u}ben) AMG has become a workhorse algorithm in scientific computing, particularly due to high-quality implementations available in standard packages~\cite{VEHenson_UMYang_2002a, RDFalgout_UMYang_2002a, BeOlSc2022}.  While it provides an efficient and robust solution algorithm for a wide class of diffusion problems, its reliance on heuristics and algorithmic parameters (that can be difficult to tune) is often seen as a difficult hurdle to overcome, particularly for critical applications.  As a result, significant effort has been invested in recent years in the development of AMG approaches with rigorous convergence bounds.  Within this area are several approaches based on the pairwise aggregation methodology~\cite{MR2914316, MR3242985, MR3069092, MR3510804, MR3716563} that offers guaranteed convergence for problems such as the graph Laplacian.  An alternative approach builds on the reduction-based multigrid methodology first proposed by Ries, Trottenberg, and Winter~\cite{MRies_UTrottenberg_GWinter_1983a}.  The reduction-based algebraic multigrid (AMGr) methodology introduced by MacLachlan et al.~\cite{SPMacLachlan_TAManteuffel_SFMcCormick_2006a} uses algebraic properties of the linear system to determine reduction-like grid-transfer operators and relaxation that again leads to guaranteed convergence rates.

The basic principle of reduction-based multigrid follows from classical cyclic reduction algorithms~\cite{PNSwarztrauber_1977a} that reduce the cost of the direct solution of linear systems, $A\vec{x} = \vec{b}$, by partitioning the degrees of freedom into two sets that we denote by $F$ and $C$ (in typical multigrid notation).  The key feature of cyclic reduction is that this partitioning should be done in a way so that the submatrix of $A$ over the set $F$, denoted $A_{FF}$, is a diagonal matrix.  (Equivalently, the set $F$ is an \textit{independent set} in the graph associated with the sparse matrix, $A$.)  Multigrid reduction~\cite{MRies_UTrottenberg_GWinter_1983a} and AMGr~\cite{SPMacLachlan_TAManteuffel_SFMcCormick_2006a} generalize this by allowing $A_{FF}$ to be non-diagonal ($F$ is not required to be an independent set), but only spectrally equivalent to a diagonal matrix.  A key question left unanswered in this work is how to generate such partitions.  MacLachlan and Saad~\cite{maclachlan2007greedy} show that the task of generating the largest possible $F$ set is an integer linear programming problem and, consequently, proposed a greedy algorithm for the partitioning with linear complexity.  Notably, they proved that if $A$ is symmetric, positive definite, and diagonally dominant, then there is a two-grid AMGr method with a guaranteed convergence bound.  Several theoretical and practical improvements to the theory and algorithms of MacLachlan and Saad~\cite{maclachlan2007greedy} have been proposed, with improved AMGr algorithms and convergence bounds~\cite{gossler2016amg, brannick2010adaptive} and improved coarsening algorithms~\cite{zaman2021coarse,ATag_etal_2021a}.  Nonsymmetric variants have also been considered~\cite{SMacLachlan_YSaad_2007b, manteuffel2018nonsymmetric, manteuffel2019nonsymmetric}.  Further research into multigrid reduction (but not AMGr) has also been commonplace in recent years, with the emergence of multigrid-reduction-in-time~\cite{Falgout_etal_2014} (MGRIT) methodologies for the solution of time-dependent PDEs and its application to multiphase flows in poromechanics~\cite{doi:10.1137/19M1256117,BUI2021114111}.

Despite these developments, attaining effective AMGr convergence has remained essentially limited to systems with matrices that are either close to diagonally dominant or can be reordered to be close to lower triangular.  Unfortunately, these limitations are significant, and preclude applying AMGr to many important and interesting classes of problems, including common AMG test cases, such as anisotropic diffusion equations, or problems discretized on anisotropic meshes.  Poor performance on anisotropic problems, in particular, is reported in existing AMGr results~\cite{maclachlan2007greedy, zaman2021coarse} that has only been overcome with expensive and impractical fixes, such as moving substantial numbers of points from $F$ to $C$, in order to construct suitably improved interpolation operators within the AMGr framework.  In this paper, we revisit the basic AMGr framework, with the goal of overcoming the barriers to achieving acceptable convergence for anisotropic problems.  To do this, we look for more practical algorithmic choices within an AMGr-style algorithm.  Specifically, we consider four ingredients for improving AMGr performance:
\begin{enumerate}
\item $C$-relaxation; while the development of AMGr by MacLachlan and co-authors~\cite{SPMacLachlan_TAManteuffel_SFMcCormick_2006a,maclachlan2007greedy} focused on $F$-relaxation, we show that using $C$-relaxation (as considered in other settings~\cite{brannick2010adaptive, Falgout_etal_2014}) can be particularly helpful for these problems;
\item Sparse Approximate Inverses (SPAI); while the original AMGr papers focused on approximating $A_{FF}^{-1}$ by a diagonal matrix, we consider using the SPAI algorithm~\cite{kolotilina1993factorized, benzi1996sparse, grote1997parallel} to offer better approximation of $A_{FF}^{-1}$ while still maintaining sparsity;
\item Strength of Connection; while the original AMGr algorithms do not rely on the classical AMG~\cite{ABrandt_SFMcCormick_JWRuge_1984a, JWRuge_KStuben_1987a} notion of ``strong connections'' to filter small entries from the matrix, we find that such filtering is critical to success for anisotropic problems, where large ``wrong-sign'' off-diagonal entries appear, but cannot be productively used in the relaxation or coarse-grid correction processes; and
\item Interpolation truncation; the combination of techniques described above leads to effective AMGr-style solvers for a wider range of problems, but with higher algorithmic complexity than is needed; thus, we use the technique of interpolation truncation~\cite{KStuben_2001a, HDeSterck_UMYang_JJHeys_2006a, HDeSterck_etal_2008b} to control these costs.
\end{enumerate}
Numerical results demonstrate that these techniques can be used in combination to improve the performance of AMGr for both isotropic and anisotropic diffusion problems.

Sparse approximate inverses are one of a class of algorithms that define preconditioners for $A\vec{x}=\vec{b}$ by prescribing a form for matrix $M$, then minimizing some norm of $I-MA$ (or $I-AM$).  Originally proposed and investigated by Benson and Frederickson~\cite{frederickson1975fast, benson1973iterative, MR683569}, recent investigations include factorized sparse approximate inverses~\cite{kolotilina1993factorized} (FSAI), which aim to compute a sparse approximation to the Cholesky factorization of SPD matrix $A$, and SPAI techniques~\cite{chow1998approximate, benzi1996sparse, grote1997parallel} that directly compute sparse approximations to $A^{-1}$.  The use of SPAI techniques in multigrid methods dates back almost to their initial introduction~\cite{frederickson1975fast, doi:10.1080/00207169008803912}, primarily to replace the use of standard relaxation schemes, such as the weighted Jacobi and Gauss-Seidel iterations.  The use of SPAI techniques for relaxation within both geometric and algebraic multigrid has been considered more recently in several ways~\cite{broker2001robust, broker2002sparse, MR1920562, broker2003parallel, wang2009multilevel, gravvanis2014algebraic, filelis2014parallel}.  Similar ideas have been used in other contexts, to build interpolation operators~\cite{wagner2000algebraic, nagel2008filtering} or improve coarse-grid operators~\cite{bolten2016sparse}.  The work in this paper is closest to the ideas presented by Bollh\"{o}fer~\cite{bollhofer2002adapted}, where SPAI was used to determine both the relaxation scheme and the interpolation operator, although the remaining details of the scheme are quite different.  We also note similarity to the work of Meurant~\cite{meurant2001numerical, meurant2002multilevel}, where entries from the AINV~\cite{benzi2000robust} preconditioner were directly used to determine interpolation alongside AINV for relaxation.

The remainder of this paper is organized as follows.  In~\cref{sec:algeb-mg-gamgr}, we give an introduction to algebraic multigrid, with a particular focus on reduction-based AMG (AMGr).  We highlight, in~\cref{sec:fail-amgr-aniso-gamgr}, that AMGr as it exists has significant difficulties for anisotropic diffusion equations, motivating the work that follows.  A key component of the algorithms considered here is the Sparse Approximate Inverse methodology of Grote and Huckle~\cite{grote1997parallel}, we review this as well in~\cref{sec:spai-gamgr}.  The main contribution of this paper is the generalized AMGr algorithm developed in~\cref{sec:gen-amgr-gamgr}.  Supporting numerical results are presented in~\cref{sec:results-gamgr}, followed by conclusions in~\cref{sec:conclusion-gamgr}.

\section{Reduction-based algebraic multigrid (AMG\lowercase{r})}\label{sec:algeb-mg-gamgr}

Multigrid methods are based on the principle of complementarity, using fine-grid relaxation and coarse-grid correction to efficiently damp all errors in the approximation of solutions to linear systems $A\vec{x}=\vec{b}$.  Geometric multigrid methods (GMG) fix a multigrid hierarchy by directly discretizing the PDE on a series of meshes defined by the problem geometry, and by adapting the relaxation scheme to complement the coarse-grid correction process defined in this way.  Algebraic multigrid methods, in contrast, do not rely on explicit knowledge of the geometry nor the PDE, instead determining the coarse levels of the multigrid hierarchy in a setup phase that precedes the solution phase of the multigrid algorithm.  In the setup, the set of degrees of freedom (or points) on the finest grid, $\Omega$, is partitioned into disjoint sets, $\Omega = C \cup F$ (with $C \cap F = \emptyset$).  The degrees of freedom in the set $C$ constitute the points on the second level.  Along with this partitioning, an interpolation operator, $P$, is constructed to map vectors from $C$ onto $\Omega$; similarly a restriction operator is defined, $R = P^T$ (in the symmetric case, as considered here), to map vectors from $\Omega$ onto $C$.  With this, the Galerkin coarse-grid operator, $A_{\ssc} = P^T A P$, is formed and the process continues recursively on $A_{\ssc}$ and $C$.
The hierarchy is constructed
until the number of nodes on a coarse grid is
sufficiently small that direct factorization of $A_{\ssc}$ is feasible. The algorithm for a two-level setup phase
is shown in~\cref{alg:AMGSetup-gamgr}. Once the setup phase
is completed, the solution phase solves the original system of
equations using a standard multigrid cycling algorithm. A two-level algorithm is shown in~\cref{alg:AMGSol-gamgr}.  Multilevel generalizations come from recursively solving $A_{\ssc}\vec{e}_{\ssc}=\vec{r}_{\ssc}$ using the two-grid methodology, either once per level (leading to a V-cycle) or multiple two-grid sweeps per level (leading, for example, to the W-cycle).  AMG methods are generally distinguished by how they define $C$ from $\Omega$, and how they define $P$ from $C$ and $A$.  Below, we review the reduction-based AMG algorithm of MacLachlan et al.~\cite{SPMacLachlan_TAManteuffel_SFMcCormick_2006a}.

\begin{algorithm}[t]
\caption{AMG Setup Phase}\label{alg:AMGSetup-gamgr}
\begin{algorithmic}[1]
\Function{amg-two-level-setup}{$A$}

\State $C, F \gets $ split the degrees of freedom into coarse and fine nodes
\State $P \gets$ form interpolation operator
\State $A_{\rm c} \gets P^T A P$
\State \textbf{return} $(P, A_{\rm c})$
\EndFunction%
\end{algorithmic}
\end{algorithm}
\begin{algorithm}[t]
\caption{AMG Solution Phase}\label{alg:AMGSol-gamgr}
\begin{algorithmic}[1]
\Function{amg-two-level-v-cycle}{$A, \vec{b}, \vec{x}, P$}

\For {$j \gets 1,\ldots,\nu_1$}\Comment{Run $\nu_1$ sweeps of pre-relaxation}
   \State $\vec{x} \gets$ relax on $\vec{x}$
\EndFor%
\State $\vec{r}_{\ssc} \gets P^{T}\left(\vec{b}-A\vec{x}\right)$
\State $\vec{e}_{\ssc} \gets$ solution of $A_{\ssc}\vec{e}_{\ssc}=\vec{r}_{\ssc}$\Comment{Solve the coarse-level problem using direct solve}
\State $\vec{x} \gets \vec{x}+P\vec{e}_{\ssc}$
\For {$j \gets 1,\ldots,\nu_2$}\Comment{Run $\nu_2$ sweeps of post-relaxation}
    \State $\vec{x} \gets$ relax on $\vec{x}$
\EndFor%
\State \textbf{return} $\vec{x}$
\EndFunction%
\end{algorithmic}
\end{algorithm}

\subsection{Reduction-based algebraic multigrid (AMGr)}\label{sec:amgr-gamgr}

Cyclic reduction~\cite{PNSwarztrauber_1977a} was originally proposed as a direct solver for certain linear systems that arose from finite-difference discretization of simple PDEs.  Assuming that the degrees of freedom are  already partitioned into coarse
and fine nodes, the linear system $A\vec{x} = \vec{b}$ is reordered to have $F$ degrees of freedom followed by $C$ degrees of freedom, writing
\begin{equation}\label{eq:permuted-gamgr}
A=
\begin{bmatrix}
 A_{\ssff} & -A_{\ssfc} \\
-A_{\sscf} &  A_{\sscc}
\end{bmatrix}
\quad
\vec{x} =
\begin{bmatrix}
  \vec{x}_{\ssf}\\
  \vec{x}_{\ssc}
\end{bmatrix}
\quad
\vec{b} =
\begin{bmatrix}
  \vec{b}_{\ssf}\\
  \vec{b}_{\ssc}
\end{bmatrix}.
\end{equation}
An \textit{exact}
algorithm for the solution of $A\vec{x} = \vec{b}$ in this partitioned
form is given by
\begin{enumerate}
  \item $\vec{y}_{\ssf} = A_{\ssff}^{-1}\vec{b}_\ssf$,
  \item Solve $\left(A_{\sscc} - A_{\sscf}A_{\ssff}^{-1}A_{\ssfc}\right)\vec{x}_{\ssc}
    = \vec{b}_{\ssc} + A_{\sscf}\vec{y}_{\ssf}$,
  \item $\vec{x}_{\ssf} = \vec{y}_{\ssf} + A_{\ssff}^{-1}A_{\ssfc}\vec{x}_{\ssc}$.
\end{enumerate}
This can be turned into an iterative method for solving
$A\vec{x}=\vec{b}$ in the usual way, replacing the right-hand side
vector, $\vec{b}$, by the evolving residual
and introducing approximations of $A_{\ssff}^{-1}$ in three places in the above algorithm, namely
\begin{equation}\label{eq:three_approx-gamgr}
  \tilde{A}_{\ssff}^{-1} \approx A_{\ssff}^{-1}
  \qquad
  \tilde{A}_{\ssc} \approx A_{\sscc} - A_{\sscf}A_{\ssff}^{-1}A_{\ssfc}
  \qquad
  W_{\ssfc} \approx A_{\ssff}^{-1}A_{\ssfc},
\end{equation}
leading to reduction-based multigrid~\cite{MRies_UTrottenberg_GWinter_1983a}.  In this form, we compute updates to the
current approximation, $\vec{x}^{(k)}$, as
\begin{enumerate}
  \item $\vec{x}^{(k+\sfrac{1}{2})}_{\ssf} = \vec{x}^{(k)}_{\ssf} +
    \tilde{A}_{\ssff}^{-1}\left(\vec{b}_{\ssf} - A_{\ssff}\vec{x}^{(k)}_{\ssf} + A_{\ssfc}\vec{x}^{(k)}_{\ssc}\right)$,
  \item Solve $\tilde{A}_{\ssc}\vec{y}_{\ssc}
    = \vec{b}_{\ssc} + A_{\sscf}\vec{x}^{(k+\sfrac{1}{2})}_{\ssf} - A_{\sscc}\vec{x}^{(k)}_{\ssc}$,
  \item $\vec{x}^{(k+1)}_{\ssc} = \vec{x}^{(k)}_{\ssc} + \vec{y}_{\ssc}$,
  \item $\vec{x}^{(k+1)}_{\ssf} = \vec{x}^{(k+\sfrac{1}{2})}_{\ssf} + W_{\ssfc}\vec{y}_{\ssc}$.
\end{enumerate}

Viewing this as a two-grid algorithm, we recognize the first step as a special form of relaxation, known as $F$-relaxation, where the approximation to $A_{\ssff}^{-1}$  is accomplished via a standard weighted Jacobi or Gauss-Seidel iteration.  The second step then represents a coarse-grid solve, where the residual is restricted to the coarse-grid by injection, and the correction, $\vec{y}_{\ssc}$, is computed using
an approximation, $\tilde{A}_{\ssc}$, of the true Schur complement, $A_{\sscc} -
A_{\sscf}A_{\ssff}^{-1}A_{\ssfc}$.  The final two steps represent the interpolation of the correction, writing the interpolation
operator $P = \left[\begin{smallmatrix} W_{\ssfc}
    \\ I \end{smallmatrix}\right]$.  This can be viewed as an approximation of the
\textit{ideal} interpolation operator~\cite{RDFalgout_PSVassilevski_2004a}, $W_{\ssfc} \approx A_{\ssff}^{-1}A_{\ssfc}$.  Notably, this algorithm differs from standard multigrid cycling in several ways, including the fixed use of injection for the restriction of the residual to the coarse grid, and the lack of post-relaxation sweeps.

As written above, there is little guidance in how to choose the three approximations in~\eqref{eq:three_approx-gamgr}.
MacLachlan et al.~\cite{SPMacLachlan_TAManteuffel_SFMcCormick_2006a} address this in their development of the reduction-based AMG (AMGr) algorithm, connecting convergence of the two-grid scheme with properties of $A_{\ssff}$.  In particular, it is assumed that $A_{\ssff}^{-1}$ can be approximated by known matrix $M_{\ssff}$ for which computing
the action of $M_{\ssff}$ on a vector is
computationally feasible. To make this rigorous, they assume that $A_{\ssff}$ can be decomposed as $A_{\ssff}=M_{\ssff}^{-1}+\mathcal{E}$, with $M_{\ssff}$ symmetric and
$0\leq\vec{x}^T\mathcal{E}\vec{x} \leq \epsilon \vec{x}^TM_{\ssff}^{-1}\vec{x}$ for all $\vec{x}$ for some $\epsilon \geq 0$, and then show that the two-grid cycle in~\cref{alg:AMGSol-gamgr} with $P=\left[
\begin{smallmatrix}
  M_{\ssff}A_{\ssfc} \\
I
\end{smallmatrix}\right]$, $A_{\ssc} = P^TAP$, and $\nu = \nu_1 = \nu_2$ pre- and post-$F$-relaxation sweeps using $M_{\ssff}$ to approximate $A_{\ssff}^{-1}$ has an error-propagation operator with norm bounded less than $1$, depending only on $\epsilon$ and $\nu$.  A technical requirement of this result (that is important below) is that the matrix $\left[ \begin{smallmatrix}
    M_{\ssff}^{-1} & -A_{\ssfc} \\
    -A_{\ssfc}^T & A_{\sscc}
  \end{smallmatrix}\right]$ must be symmetric and positive semi-definite for the convergence result to hold.

While this work is insightful, it does not address the fundamental question of how to generate a partitioning for which the assumptions hold with small parameter $\epsilon$.  To answer this question, MacLachlan and Saad~\cite{maclachlan2007greedy}
propose to partition the rows and columns of $A$
in order to ensure the diagonal dominance
of $A_{\ssff}$, allowing $M_{\ssff}$ to be chosen as a diagonal matrix. In
particular, for each row, $i$, the diagonal dominance of row $i$ over the $F$ points is quantified by
\[
  \eta_i=\frac{|A_{ii}|}{\sum_{j\in F}|A_{ij}|}.
\]
Then, $A_{\ssff}$ is said to be $\eta$-diagonally dominant if $\eta_i \geq \eta$
for all $i \in F$, for some $\eta>1/2$ that measures the diagonal dominance
of $A_{\ssff}$.  If $A_{\ssff}$ is $\eta$-diagonally
dominant, then the diagonal matrix, $M_{\ssff}$, with
$(M_{\ssff})_{ii}^{-1} = (2-\frac{1}{\eta})A_{ii}$ for all $i\in F$ yields $0\leq \vec{x}^T\mathcal{E}\vec{x} \leq
\frac{2-2\eta}{2\eta-1}\vec{x}^TM_{\ssff}^{-1}\vec{x}$ for all $\vec{x}$,
giving $\epsilon = \frac{2-2\eta}{2\eta-1}$.  Furthermore, if $A$ is symmetric, positive-definite, and
diagonally dominant, then this condition guarantees
that all conditions of the theory from MacLachlan et al.~\cite{SPMacLachlan_TAManteuffel_SFMcCormick_2006a} are satisfied.

In
addition to establishing this connection between the diagonal dominance
parameter $\eta$ and the convergence parameter, $\epsilon$,
MacLachlan and Saad~\cite{maclachlan2007greedy} consider practical algorithms for computing a partition with $\eta$-diagonally dominant $A_{\ssff}$.  They show that finding the largest $F$-set with $\eta$-diagonally dominant $A_{\ssff}$ is an NP-complete problem, but propose a
greedy algorithm to approximately solve the optimization problem for the largest such $F$-set.  The greedy algorithm acts iteratively, adding points to the $C$-set one at a time, and moving any points that are guaranteed to satisfy the diagonal dominance constraint into the $F$-set, until a full partition is computed.
While the greedy coarsening algorithm was demonstrated to be effective in some settings, Zaman et al.~\cite{zaman2021coarse} demonstrate that there are also cases where the resulting optimality gap can be significant.  To address this, they propose to apply simulated annealing to the same optimization problem, showing that this approach can produce substantially better partitionings than the greedy approach, albeit at the greatly increased cost of many simulated annealing steps.  Similar work by Taghibakhshi et al.~\cite{ATag_etal_2021a} uses reinforcement learning to solve the same problem at a lower cost.

Several generalizations of both the theory and practice of reduction-based multigrid methods have also been developed.  A generalization to non-symmetric M-matrices was proposed and analyzed by Mense and Nabben~\cite{mense2008algebraic}, using the tools of weak regular splitting~\cite{RSVarga_2000}.  For symmetric and positive definite problems, Brannick et al.~\cite{brannick2010adaptive} study the introduction of more general relaxation schemes, as well as the use of different approximations of $A_{\ssff}$ for interpolation and relaxation.  Gossler and Nabben~\cite{gossler2016amg} examine generalization of AMGr to the use of Chebyshev polynomial acceleration of multiple relaxation sweeps.  For strongly non-symmetric systems, Manteuffel et al.~\cite{manteuffel2018nonsymmetric,manteuffel2019nonsymmetric} have proposed similar approaches using so-called approximate ideal restriction (AIR) techniques, that offer excellent performance for advection-dominated problems.  None of the above schemes, however, address the poor performance observed in AMGr-type methods for anisotropic problems, which is the motivation for the present work. These advances, coupled with the success of other approaches that make use of multigrid reduction principles (but not, specifically, AMGr)~\cite{Falgout_etal_2014,doi:10.1137/19M1256117,BUI2021114111}, suggest that further investigation of the AMGr paradigm may be worthwhile.  Here, in particular, we find that AMGr robustness and performance can be greatly enhanced by conforming less strictly to the theoretical framework of MacLachlan et al.~\cite{SPMacLachlan_TAManteuffel_SFMcCormick_2006a} and making use of some techniques from classical AMG.

\subsection{Failure of AMGr for anisotropic diffusion}\label{sec:fail-amgr-aniso-gamgr}

To demonstrate the convergence problems, we consider applying AMGr (following the prescription of MacLachlan and Saad~\cite{maclachlan2007greedy}) to the solution of the two-dimensional anisotropic diffusion problem,
\begin{equation}\label{eq:anisotropic_diffusion-gamgr}
-\nabla \cdot \boldsymbol{K}(x, y)\nabla u(x, y) = b(x,y)
\end{equation}
in the
domain $[0,1]\times [0,1]$ with Dirichlet boundary conditions.
We choose the tensor coefficient $\boldsymbol{K}(x, y)=QHQ^{T}$,
where
\begin{equation}\label{eq:Q-H-def-gamgr}
  Q=\begin{bmatrix}
\cos(\theta) & -\sin(\theta)
\\ \sin(\theta) & \cos(\theta)
 \end{bmatrix},
 \qquad
H=\begin{bmatrix}
\delta_1 & 0
\\ 0 & \delta_2
 \end{bmatrix},
\end{equation}
where $\theta$ specifies the
direction of anisotropy in the problem, and $\delta_1$ and $\delta_2$ specify the strengths. We consider $\delta_1=10^{-6}$ and $\delta_2=1$ for this problem. For $\theta = 0$ this gives the grid-aligned anisotropic equation $-10^{-6} u_{xx} - u_{yy} = b$, while $0 < \theta < \pi/2$ gives a non-grid-aligned diffusion tensor.  \Cref{tab:conv_anisot_struc_geom-lloyd_2lev-gamgr} presents convergence results for the standard finite-difference and bilinear finite element discretizations of this problem, with $\theta = 0$, $\pi/6$, and $\pi/4$.  For simplicity, we present results for a uniform $32\times 32$ grid, although similar results are observed for larger meshes.  Here, and in all results that follow, we measure convergence by solving the homogeneous problem, $A\vec{x}=\vec{0}$, with a randomly chosen initial guess for $\vec{x}$.  Writing $\vec{e}^{(k)}$ as the error in the $k^{\text{th}}$ approximation to $\vec{x} = \vec{0}$, we estimate the asymptotic convergence factor by running 50 (stationary) multigrid iterations, then estimating $\rho \approx \left(\|\vec{e}^{(50)}\|_A / \|\vec{e}^{(10)}\|_A\right)^{1/40}$, averaging convergence over the final 40 iterations.  We note that, in all cases, the coarsening algorithm (simulated annealing, in this case) generates a partitioning such that the matrix $A_{\ssff}$ is well-approximated by a diagonal matrix, $M_{\ssff}$; however, only in the case of the finite-difference discretization of the grid-aligned diffusion equation (when the discretization matrix, $A$, \textit{is} diagonally dominant), does the required semidefiniteness of $\left[ \begin{smallmatrix}
    M_{\ssff}^{-1} & -A_{\ssfc} \\
    -A_{\ssfc}^T & A_{\sscc}
\end{smallmatrix}\right]$ hold.  This correlates strongly with the resulting measured asymptotic convergence factor for the method.  Both MacLachlan and Saad~\cite{maclachlan2007greedy} and Zaman et al.~\cite{zaman2021coarse} consider remedies for this behavior, such as augmenting the $C$ set in a style similar to classical AMG\@; while this improves the overall convergence of the method, it also leads to greatly increased grid and operator complexities, making it an unsatisfactory solution.

\begin{table}
	\caption{Performance of two-level AMGr for the anisotropic diffusion problem.}
    \centering
		\begin{tabular}{c  c  c  c  c  N  \Cgrid\CopTwo}
			\toprule
      \multirow{2}{*}{$\theta$} &
      \multirow{2}{*}{Discretization} &
      \multicolumn{2}{c}{Eigenvalues of $M_{\ssff}A_{\ssff}$} &
      \multirow{2}{*}{$\left[ \begin{smallmatrix} M_{\ssff}^{-1} & -A_{\ssfc} \\ -A_{\ssfc}^T & A_{\sscc} \end{smallmatrix}\right]$} &
      Convergence factor &
      \multicolumn{2}{c}{Complexities} \\
       & & min & max &  & $\rho$ & $C_{\text{grid}}$ & $C_{\text{op}}$ \EndTableHeader\\
      \midrule
                    \multirow{2}{*}{0} & FD & 1.00 & 3.00 & positive definite & 0.75 & 1.33 & 1.49\\
                    & FE & 1.00 & 6.83 & indefinite & 0.98 & 1.49 & 2.14\\
                    \cmidrule[0.25pt](lr){1-8}
                    \multirow{2}{*}{$\pi/6$} & FD & 1.04 & 5.59 & indefinite & 0.96 & 1.42 & 1.87\\
                     & FE & 1.00 & 7.09 & indefinite & 0.97 & 1.36 & 1.72 \\
                    \cmidrule[0.25pt](lr){1-8}
                    \multirow{2}{*}{$\pi/4$} & FD & 1.21 & 6.73 & indefinite & 0.96 & 1.39 & 1.79 \\
                     & FE & 1.12 & 5.80 & indefinite & 0.96 & 1.33 & 1.59 \\
      \bottomrule
		\end{tabular}\label{tab:conv_anisot_struc_geom-lloyd_2lev-gamgr}
\end{table}

Here, and in all tables that follow, we use color-coding to indicate quality of the results shown.  For measured convergence factors, we denote a ``good'' convergence factor to be below 0.4 (indicated in green), while a ``bad'' convergence factor, above 0.8, is shown in red (with values in between, $0.4 < \rho < 0.8$, shown in black text). An important consideration in assessing the quality of coarsening is the resulting complexity of the multigrid algorithm, as this will vary with the coarse grids chosen in algebraic multigrid.  We use two common measures: grid and operator complexity.  The AMG grid complexity, $C_{\rm{grid}}$, is the ratio of the sum of the number of DoFs on each level of multigrid hierarchy (including the finest) to that on the finest level. Similarly, the operator complexity, $C_{\rm{op}}$, is the ratio of the sum of the number of nonzeros in the system matrices on each level of hierarchy (including the finest) to that on the finest level. For two-grid operator complexity, we highlight results in green if the complexity is below 1.5, in red if it is above 2.5, and in orange for values between 2.0 and 2.5.  As operator complexity, in particular, is expected to grow with the number of levels in the hierarchy, we use similar highlighting with different thresholds for three-grid and multigrid operator complexity, showing results in green if it is below 2.0, red if it is above 3.0, and orange for values between 2.5 and 3.0.  As the results that follow show relatively little variation in AMG grid complexity, we choose not to highlight values for this measure.

\section{Sparse Approximate Inverse (SPAI) methods}\label{sec:spai-gamgr}

While originally proposed in the 1970's by Benson and Frederickson~\cite{frederickson1975fast, benson1973iterative}, SPAI techniques were more systematically developed and studied in the 1990's (and subsequently) by a number of authors~\cite{kolotilina1993factorized, benzi1996sparse, grote1997parallel, chow1998approximate, benzi1999comparative}.  The general idea of SPAI techniques is to compute a matrix, $M$, to minimize some norm of $I-MA$ or $I-AM$, with constraints on the sparsity of $M$.  These constraints may be fixed (e.g., some fixed set of elements of $M$ is allowed to be nonzero), or may be adaptively determined by trying to best minimize the chosen norm within some limitations on either the total number of nonzero entries in $M$ or the row/column-wise number of nonzero elements.  Here, we focus on the variant of the SPAI algorithm proposed by Hawkins and Chen~\cite{doi:10.1137/S1064827503423500}, in which the Frobenius norm of $B-AM$ is minimized for a given matrix, $B$, over a fixed nonzero pattern for each column.  If $B=I$, then this reduces to a simplified version of the SPAI algorithm of Grote and Huckle~\cite{grote1997parallel}, omitting their adaptive calculation for increasing the nonzero pattern for each column.

\Cref{alg:SPAI-gamgr} presents the SPAI algorithm, where the inputs are given by matrices $A$ and $B$,
and a nonzero sparsity pattern, $\mathcal{S}$, for the sparse approximate inverse $M$.
The algorithm loops independently over each column, $j$, in $M$.  In the initialization stage of the algorithm (Lines~\ref{line:start_init} through~\ref{line:end_init}), the rows, $\mathcal{J}$, in the initial sparsity pattern of $\mathcal{S}$ for column $j$ are extracted, as is the set of rows, $\mathcal{I}$, of $A$ for which matrix $A$ has a nonzero entry in a column in $\mathcal{J}$, defining a submatrix, $\bar{A}$, of $A$ that is used to initialize column $j$ of $M$.  Two auxiliary vectors are also formed, corresponding to the full $j^{\text{th}}$ column of $B$, denoted $\vec{b}$, and its restriction to the rows of $\bar{A}$, denoted $\bar{\vec{b}}$.  Column $j$ of $M$ then comes from using the QR decomposition of $\bar{A}$ to solve the unconstrained minimization problem of minimizing $\left\|\bar{\vec{b}}-\bar{A}\bar{\vec{m}}\right\|_2$, noting that $\bar{A}$ is expected, by its construction, to have more rows than columns, so that this is not expected to yield a zero residual.  The computed solution, $\bar{\vec{m}}$, is injected into a full vector, $\vec{m}$, which becomes the $j^{\text{th}}$ column of $M$.

\begin{algorithm}
  \caption{SPAI algorithm}\label{alg:SPAI-gamgr}
  \begin{algorithmic}[1]
  \Function{spai}{$A, B, \mathcal{S}$}
    \For {$j \gets 1,\ldots,n_{\text{col}}$}   \Comment{$n_{\text{col}}$ is number of columns in $A$}
        \State $\mathcal{J} \gets \{i \mid (i,j) \in \mathcal{S}\}$ \label{line:start_init}
        \State $\mathcal{I} \gets$ set of indices of nonzero rows of $A(:,\mathcal{J})$
        \State $\vec{b} \gets B(:, j)$
        \State $\bar{A} \gets A(\mathcal{I}, \mathcal{J})$
        \State $\bar{\vec{b}} \gets B(\mathcal{I}, j)$\label{line:end_init}
        \State Compute QR decomposition of $\bar{A}$
        \State $\bar{\vec{m}} \gets \argmin_{\bar{\vec{m}}} \left\lVert\bar{\vec{b}}-\bar{A}\bar{\vec{m}}\right\rVert_2$ \Comment{Unconstrained least-squares via QR}
        \State $\vec{m} \gets \bar{\vec{m}}$ with inserted zeros
        \State $M(:,j) \gets \vec{m}$
    \EndFor%
    \State \textbf{return} $M$
    \EndFunction%
  \end{algorithmic}
\end{algorithm}

Sparse approximate inverse algorithms similar to~\cref{alg:SPAI-gamgr} have been investigated for use in both relaxation and interpolation in several settings in the past.  However, this usage has generally been in defining approximations to the inverse of $A$ in its entirety, while we look at the possible use of SPAI techniques through the lens of the AMGr methodology.  In what follows, we will make use of SPAI in three ways:
\begin{enumerate}
\item In $F$-relaxation where, given a proxy matrix, $\hat{A}_{\ssff}$, for $A_{\ssff}$ (possibly equal to $A_{\ssff}$), and $B = I_{\ssff}$, $M_{\ssff}$ is constructed as the SPAI approximation to $\hat{A}_{\ssff}^{-1}$ with a fixed sparsity pattern equal to that of $\hat{A}_{\ssff}$;
\item In $C$-relaxation where, given a proxy matrix, $\hat{A}_{\sscc}$, for $A_{\sscc}$ (possibly equal to $A_{\sscc}$), and $B = I_{\sscc}$, $M_{\sscc}$ is constructed as the SPAI approximation to $\hat{A}_{\sscc}^{-1}$ with a fixed sparsity pattern equal to that of $\hat{A}_{\sscc}$; and
\item In interpolation, where we use the Hawkins and Chen modfication~\cite{doi:10.1137/S1064827503423500}, to solve $\min_{W} \|\hat{A}_{\ssfc} - \hat{A}_{\ssff}W\|_F$ for sparse approximate columns of $W = \hat{A}_{\ssff}^{-1}\hat{A}_{\ssfc}$ for proxy matrices, $\hat{A}_{\ssff}$ and $\hat{A}_{\ssfc}$, for $A_{\ssff}$ and $A_{\ssfc}$, respectively, with a fixed nonzero pattern equal to that of $\hat{A}_{\ssfc} + \hat{A}_{\ssff}\hat{A}_{\ssfc}$.
\end{enumerate}
The first two of these can be viewed as generalizations of the use of SPAI on all of $A$ as relaxation~\cite{broker2002sparse, broker2003parallel} to the $F$- and $C$-relaxations typically used in reduction-based AMG\@.  The third bears similarity to Meurant's Algorithm I3~\cite{meurant2001numerical}, where SPAI on $A$ (either in its original ordering or reordered according to the $F$-$C$ partitioning) is used to generate an approximate inverse matrix, $M$, from which $M_{\ssff}$ is extracted to form an interpolation operator $\left[ \begin{smallmatrix}
    M_{\ssff}A_{\ssfc} \\
    I
\end{smallmatrix}\right]$.  We note that this is akin to approximating ideal interpolation, $\left[\begin{smallmatrix} A_{\ssff}^{-1}A_{\ssfc} \\ I \end{smallmatrix} \right]$ by using an approximation to $\left(A^{-1}\right)_{\ssff}$ rather than $\left(A_{\ssff}\right)^{-1}$.  As discussed below, direct use of SPAI to approximate $A_{\ssff}^{-1}$, $A_{\sscc}^{-1}$, and $A_{\ssff}^{-1}A_{\ssfc}$ in these contexts does not directly lead to effective performance, so we consider additional tools from standard AMG development to both improve convergence and lower cost.

\section{Generalizing AMGr}\label{sec:gen-amgr-gamgr}

{\bf Baseline results:} From the results in~\cref{tab:conv_anisot_struc_geom-lloyd_2lev-gamgr} and those documented in other works~\cite{maclachlan2007greedy, zaman2021coarse, ATag_etal_2021a}, the coarse-grid correction process emerges as a primary source for the poor performance of AMGr on anisotropic problems.  Using $M_{\ssff}A_{\ssfc}$ with diagonal $M_{\ssff}$ for the $C$-to-$F$ interpolation matrix is more restrictive than the interpolation operators used in classical multigrid, as interpolation to an $F$ point is only allowed from directly connected $C$ points (corresponding to the nonzero entries in $A_{\ssfc}$).  To test this theory, we first use SPAI to determine an interpolation operator of the form $P = \left[ \begin{smallmatrix}
	W \\ I
  \end{smallmatrix} \right]$, where the sparsity pattern of $W$ is fixed to match that of $A_{\ssfc} + A_{\ssff}A_{\ssfc}$, allowing for interpolation to a fine-grid point from both directly connected $C$ points and $C$ points that are directly connected to an adjacent $F$ point.  At the same time, we replace relaxation based on a diagonal stencil with weighted relaxation using the SPAI approximation to $A_{\ssff}^{-1}$, denoted $M_{\ssff}$, with the sparsity pattern of $A_{\ssff}$, and weight $\sigma_{\ssf} = 2/(\lambda_{\text{min}}+\lambda_{\text{max}})$ for extremal eigenvalues, $\lambda_{\text{min}}$ and $\lambda_{\text{max}}$ of $M_{\ssff}A_{\ssff}$, chosen to minimize the spectral radius of $I-\sigma_{\ssf}M_{\ssff}A_{\ssff}$. (Preliminary results (not shown here), replacing only interpolation and not relaxation, show qualitatively similar results to those below in~\cref{tab:aniso_FE_approx_ideal_interp_Frelax_cb3-gamgr}, but with notably larger convergence factors for the $\theta =0$ and $\theta = \pi/6$ cases.)

To test the effects of these modifications, we again consider the anisotropic diffusion equation given in~\eqref{eq:anisotropic_diffusion-gamgr}, for three angles, $\theta = 0$, $\pi/6$, and $\pi/4$, with convergence factors shown in~\cref{tab:aniso_FE_approx_ideal_interp_Frelax_cb3-gamgr}.  To decouple the impact of these choices from that of the coarsening, we use a geometric coarse grid chosen as semi-coarsening by a factor of three in the $y$-direction (the direction of strong connections in these cases).  We observe significant improvement in convergence in both the $\theta = 0$ and $\theta = \pi/6$ cases, in comparison to the results in~\cref{tab:conv_anisot_struc_geom-lloyd_2lev-gamgr}, although performance for $\theta = \pi/4$ is much worse.  Yet, we also note that the complexity of these cycles is high, with two-grid operator complexities above 2.0 due to many small nonzero entries in the resulting interpolation operators that lead to large numbers of nonzero entries in the Galerkin coarse-grid operator, $P^T A P$. Semi-coarsening by a factor of three is used as preliminary experiments with algebraic coarsening indicated coarsenings with similar grid complexity were attained with typical parameters, as will be seen below.

\begin{table}
  \caption{Two-level AMGr convergence factors for anisotropic FE discretization using semi-coarsening in the $y$ direction by a factor of three. Interpolation of $F$-nodes
 uses the SPAI approximation to $A_{\ssff}^{-1}A_{\ssfc}$ and the SPAI approximation to $A_{\ssff}^{-1}$ is used for relaxation.}
	\centering
		\begin{tabular}{c|N\Cgrid\CopTwo|N\Cgrid\CopTwo|N\Cgrid\CopTwo}  %
			\toprule
			      & \multicolumn{3}{c|}{$\theta=0$} & \multicolumn{3}{c|}{$\theta=\pi/6$} & \multicolumn{3}{c}{$\theta=\pi/4$}\\
     Grid size & $\rho$ & $C_{\text{grid}}$ & $C_{\text{op}}$ & $\rho$ & $C_{\text{grid}}$ & $C_{\text{op}}$ & $\rho$ & $C_{\text{grid}}$ & $C_{\text{op}}$ \EndTableHeader\\
     \midrule
                        $16\times16$ & 0.021 & 1.31 & 1.90 & 0.006 & 1.31 & 1.90 & 0.213 & 1.31 & 1.90 \\
			            $32\times32$ & 0.023 & 1.31 & 2.02 & 0.019 & 1.31 & 2.02 & 0.518 & 1.31 & 2.02 \\
                        $64\times64$ & 0.023 & 1.33 & 2.14 & 0.065 & 1.33 & 2.14 & 0.808 & 1.33 & 2.14 \\
			          $128\times128$ & 0.024 & 1.33 & 2.17 & 0.210 & 1.33 & 2.17 & 0.921 & 1.33 & 2.17 \\
			\bottomrule
    \end{tabular}\label{tab:aniso_FE_approx_ideal_interp_Frelax_cb3-gamgr}
\end{table}

{\bf Strong connections:} To address the high cost encountered in the baseline, we introduce strong connections into the algorithm.  A typical row in the matrix for an anisotropic diffusion operator contains both small entries and large but ``wrong-sign'' entries, where there are positive contributions in directions other than the strong direction of diffusion in the PDE\@.  In response, we introduce a filtering stage, where we compute a \emph{proxy} matrix, $\hat{A}$, for the given system matrix, $A$.  We first compute strong connections using the classical Ruge-St\"{u}ben definition of strength of connection, defining point $i$ to be strongly connected to point $j$ if
\begin{equation}\label{eq:anisostrength-gamgr}
  -A_{ij} \geq \frac{1}{2} \max_{k\neq i} -A_{ik},
\end{equation}
where a strength parameter of 1/2 is selected as is typical for anisotropic PDEs and where only negative off-diagonal entries are allowed as strong connections (also common practice).

For anisotropic diffusion equations discretized by bilinear finite elements on uniform grids, \cref{eq:anisostrength-gamgr} results in two strong connections for each interior node, aligned vertically (north and south), for $\theta = 0$, two strong connections in the north-east and south-west directions for $\theta = \pi/4$, but four strong connections for $\theta = \pi/6$, including north, south, north-east, and south-west points.  To preserve the row-sum that is typically needed for best AMG performance, we define $\hat{A}$ to have off-diagonal entries matching those of $A$ for strong connections, and diagonal entries adjusted by subtracting any weak connections in each row of $A$ from its diagonal value (so-called ``lumping'' of the weak connections to the diagonal, as in classical AMG).

Next, we repeat the experiments above, but define interpolation as the SPAI approximation to $\hat{A}_{\ssff}^{-1}\hat{A}_{\ssfc}$ and use $F$-relaxation based on the SPAI approximation to $\hat{A}_{\ssff}^{-1}$, with results presented in~\cref{tab:aniso_FE_newAMGr_no_crelax_no_scaling_cb3-gamgr}.  For $\theta = 0$ and $\theta = \pi/4$, each $F$-point has a single strongly-connected $C$ neighbor and a single strongly-connected $F$ neighbor, while each $F$-point for $\theta=\pi/6$ has two of each.
Using semi-coarsening in the $y$-direction by a factor of three, this results in interpolation to each $F$-point from two $C$-points for $\theta = 0$ and $\pi/4$ and from five $C$-points for $\theta = \pi/6$ (where the two strongly connected $F$ neighbors of an $F$-point have a total of three strongly connected $C$ neighbors).  From the table, we see that using this definition of strength results in improved and grid-independent convergence for the case of $\theta = \pi/4$, and degraded (but still grid-independent) convergence for $\theta = 0$.  However, significant degradation in convergence occurs for the case of $\theta = \pi/6$.  Nonetheless, the use of strong connections has greatly improved the two-grid operator complexities, particularly for the $\theta=0$ case, where it now matches that of geometric multigrid.
\begin{table}
  \caption{Two-level AMGr convergence factors for anisotropic FE discretization using semi-coarsening in the $y$ direction by a factor of three.
    Interpolation of $F$-nodes is based on
a SPAI approximation to $\hat{A}_{\ssff}^{-1}\hat{A}_{\ssfc}$ and relaxation uses a SPAI approximation to $\hat{A}_{\ssff}^{-1}$, {\bf where $\hat{A}$ is the ``lumped'' matrix of strong connections computed from $A$}.}
	\centering
		\begin{tabular}{c|N\Cgrid\CopTwo|N\Cgrid\CopTwo|N\Cgrid\CopTwo}  %
			\toprule
			      & \multicolumn{3}{c|}{$\theta=0$} & \multicolumn{3}{c|}{$\theta=\pi/6$} & \multicolumn{3}{c}{$\theta=\pi/4$}\\
     Grid size & $\rho$ & $C_{\text{grid}}$ & $C_{\text{op}}$ & $\rho$ & $C_{\text{grid}}$ & $C_{\text{op}}$ & $\rho$ & $C_{\text{grid}}$ & $C_{\text{op}}$ \EndTableHeader\\
     \midrule
                        $16\times16$ & 0.351 & 1.31 & 1.28 & 0.289 & 1.31 & 1.60 & 0.696 & 1.31 & 1.42 \\
			            $32\times32$ & 0.359 & 1.31 & 1.30 & 0.487 & 1.31 & 1.66 & 0.718 & 1.31 & 1.47 \\
                        $64\times64$ & 0.359 & 1.33 & 1.32 & 0.745 & 1.33 & 1.73 & 0.716 & 1.33 & 1.52 \\
			          $128\times128$ & 0.359 & 1.33 & 1.32 & 0.906 & 1.33 & 1.75 & 0.716 & 1.33 & 1.53 \\
			\bottomrule
    \end{tabular}\label{tab:aniso_FE_newAMGr_no_crelax_no_scaling_cb3-gamgr}
\end{table}

{\bf Interpolation Scaling:} \Cref{tab:aniso_FE_approx_ideal_interp_Frelax_cb3-gamgr,tab:aniso_FE_newAMGr_no_crelax_no_scaling_cb3-gamgr} underscore the potential of SPAI-based AMGr for anisotropic diffusion equations, but also highlight that acceptable and scalable convergence is not robust.
From the poor convergence for $\theta = \pi/6$ in~\cref{tab:aniso_FE_newAMGr_no_crelax_no_scaling_cb3-gamgr}, we found that even
if the ``lumped'' matrix, $\hat{A}$, used to form interpolation retains the
property that rows away from boundary conditions have zero row sum (and that
$\hat{A}$ is an M-matrix), the interpolation operator determined by SPAI does
not accurately interpolate the coarse-grid constant function onto the fine-grid.
This is not surprising, since SPAI computes the
interpolation operator column-wise, yet interpolation to any fixed fine-grid
vector is a row-wise property of matrix $P$; however, it does indicate a potential reason for the degraded convergence observed in~\cref{tab:aniso_FE_newAMGr_no_crelax_no_scaling_cb3-gamgr}, as classical AMG is well-known to be ineffective when global near null-space modes, such as the constant function, are not well-approximated by the range of interpolation.

To address the lack of constant interpolation,
we post-process the interpolation
generated by SPAI, using left diagonal scaling of $W\approx
\hat{A}_{\ssff}^{-1}\hat{A}_{\ssfc}$ so that $\vec{1}_{\ssf} =
SW\vec{1}_{\ssc}$.  This is accomplished by computing $\vec{s} =
W\vec{1}_{\ssc}$, followed by defining diagonal matrix $S$ with entries $1/s_i$ on its
diagonal.  This yields the convergence factors in~\cref{tab:aniso_FE_newAMGr_no_crelax_const_scaling_cb3-gamgr}.
The results offer concrete improvement over~\cref{tab:aniso_FE_approx_ideal_interp_Frelax_cb3-gamgr,tab:aniso_FE_newAMGr_no_crelax_no_scaling_cb3-gamgr}, in that they offer
scalable convergence for all three problems (without increasing grid or operator complexities).  Even so, the overall
convergence factors between 0.65 and 0.8 are insufficient to be considered
an effective AMG solver.

One possible cause of the degraded convergence
is poor interpolation near Dirichlet boundaries, where the
constant vector is not an accurate indicator of the slowest-to-converge modes of
relaxation (or ``algebraically smooth errors'' in classical AMG).  As a remedy,
we use a similar scaling of
interpolation, that we refer to as ``improved iteration'' scaling, running a
set number of sweeps of (full grid) weighted Jacobi relaxation on the
homogeneous problem with the constant vector as an initial guess, to produce a
relaxed vector, $\vec{z}$, and followed by a similar diagonal scaling computed to
ensure that $\vec{z}_{\ssf} = SW\vec{z}_{\ssc}$.  With this modification and five sweeps of relaxation,~\cref{tab:aniso_FE_newAMGr_no_crelax_improved_iter_scaling_cb3-gamgr}
shows notable improvement in
performance for both the case of $\theta = 0$ and $\theta = \pi/6$, but still
disappointing (albeit grid-independent) convergence for $\theta = \pi/4$.

\begin{table}
  \caption{Two-level AMGr convergence factors for anisotropic FE discretization using semi-coarsening in the $y$ direction by a factor of 3. Here, $F$-nodes are
  interpolated using the SPAI approximation to $\hat{A}_{\ssff}^{-1}\hat{A}_{\ssfc}$, {\bf postprocessed to exactly interpolate the constant vector}, and the SPAI approximation to $\hat{A}_{\ssff}^{-1}$ is used for relaxation, where $\hat{A}$ is the ``lumped'' matrix of strong connections computed from $A$.}
	\centering
		\begin{tabular}{c|N\Cgrid\CopTwo|N\Cgrid\CopTwo|N\Cgrid\CopTwo}  %
			\toprule
			      & \multicolumn{3}{c|}{$\theta=0$} & \multicolumn{3}{c|}{$\theta=\pi/6$} & \multicolumn{3}{c}{$\theta=\pi/4$}\\
     Grid size & $\rho$ & $C_{\text{grid}}$ & $C_{\text{op}}$ & $\rho$ & $C_{\text{grid}}$ & $C_{\text{op}}$ & $\rho$ & $C_{\text{grid}}$ & $C_{\text{op}}$ \EndTableHeader\\
     \midrule
                        $16\times16$ & 0.751 & 1.31 & 1.28 & 0.641 & 1.31 & 1.60 & 0.765 & 1.31 & 1.42 \\
			            $32\times32$ & 0.776 & 1.31 & 1.30 & 0.640 & 1.31 & 1.66 & 0.751 & 1.31 & 1.47 \\
                        $64\times64$ & 0.772 & 1.33 & 1.32 & 0.649 & 1.33 & 1.73 & 0.744 & 1.33 & 1.52 \\
			          $128\times128$ & 0.772 & 1.33 & 1.32 & 0.656 & 1.33 & 1.75 & 0.741 & 1.33 & 1.53 \\
			\bottomrule
    \end{tabular}\label{tab:aniso_FE_newAMGr_no_crelax_const_scaling_cb3-gamgr}
\end{table}
\begin{table}
  \caption{Two-level AMGr convergence factors for anisotropic FE discretization using semi-coarsening in the $y$ direction by a factor of 3. Here, $F$-nodes are interpolated using the SPAI approximation to $\hat{A}_{\ssff}^{-1}\hat{A}_{\ssfc}$, {\bf postprocessed to exactly interpolate a relaxed vector}, and the SPAI approximation to $\hat{A}_{\ssff}^{-1}$ is used for relaxation, where $\hat{A}$ is the ``lumped'' matrix of strong connections computed from $A$.}
	\centering
		\begin{tabular}{c|N\Cgrid\CopTwo|N\Cgrid\CopTwo|N\Cgrid\CopTwo}  %
			\toprule
			      & \multicolumn{3}{c|}{$\theta=0$} & \multicolumn{3}{c|}{$\theta=\pi/6$} & \multicolumn{3}{c}{$\theta=\pi/4$}\\
     Grid size & $\rho$ & $C_{\text{grid}}$ & $C_{\text{op}}$ & $\rho$ & $C_{\text{grid}}$ & $C_{\text{op}}$ & $\rho$ & $C_{\text{grid}}$ & $C_{\text{op}}$ \EndTableHeader\\
     \midrule
                        $16\times16$ & 0.379 & 1.31 & 1.28 & 0.177 & 1.31 & 1.60 & 0.720 & 1.31 & 1.42 \\
			            $32\times32$ & 0.382 & 1.31 & 1.30 & 0.197 & 1.31 & 1.66 & 0.719 & 1.31 & 1.47 \\
                        $64\times64$ & 0.379 & 1.33 & 1.32 & 0.256 & 1.33 & 1.73 & 0.717 & 1.33 & 1.52 \\
			          $128\times128$ & 0.377 & 1.33 & 1.32 & 0.284 & 1.33 & 1.75 & 0.717 & 1.33 & 1.53 \\
			\bottomrule
    \end{tabular}\label{tab:aniso_FE_newAMGr_no_crelax_improved_iter_scaling_cb3-gamgr}
\end{table}

{\bf $C$-relaxation:} As a final modification we employ the use of $C$-relaxation alongside $F$-relaxation.  As has been considered in MGRIT~\cite{Falgout_etal_2014} and other contexts, replacing simple $F$-relaxation with sweeps of $FCF$-relaxation (that is, relaxation over the $F$-points, followed by relaxation over the $C$-points, then again over the $F$-points, with updated residual values between each sweep) is known to greatly improve multigrid performance in some settings.  Results using $FCF$-relaxation are shown in~\cref{tab:aniso_FE_newAMGr_improved_iter_scaling_crelax_cb3-gamgr}, where $C$-relaxation is again computed with SPAI on $\hat{A}_{\sscc}$.  We see that including $C$-relaxation results in a dramatic effect for $\theta = \pi/4$, reducing the convergence factor to nearly 0.1, and a notable effect for $\theta = 0$.  For $\theta = \pi/6$, adding $C$-relaxation has little influence on convergence, noting it does not harm performance. We note that another reasonable option would be to run a sweep of full-grid relaxation before the $F$-relaxation, and that this may be equivalent to $FC$-relaxation for certain relaxation schemes and partitions, but we do not explore this possibility here.
\begin{table}
  \caption{Two-level AMGr convergence factors for anisotropic FE discretization using semi-coarsening in the $y$ direction by a factor of 3. Here, $F$-nodes are
  interpolated using the SPAI approximation to $\hat{A}_{\ssff}^{-1}\hat{A}_{\ssfc}$, postprocessed to exactly interpolate a relaxed vector, and the SPAI approximations to $\hat{A}_{\ssff}^{-1}$ and $\hat{A}_{\sscc}^{-1}$ are used {\bf for $FCF$-relaxation}, where $\hat{A}$ is the ``lumped'' matrix of strong connections computed from $A$.}
	\centering
		\begin{tabular}{c|N\Cgrid\CopTwo|N\Cgrid\CopTwo|N\Cgrid\CopTwo}  %
			\toprule
			      & \multicolumn{3}{c|}{$\theta=0$} & \multicolumn{3}{c|}{$\theta=\pi/6$} & \multicolumn{3}{c}{$\theta=\pi/4$}\\
     Grid size & $\rho$ & $C_{\text{grid}}$ & $C_{\text{op}}$ & $\rho$ & $C_{\text{grid}}$ & $C_{\text{op}}$ & $\rho$ & $C_{\text{grid}}$ & $C_{\text{op}}$ \EndTableHeader\\
     \midrule
                        $16\times16$ & 0.233 & 1.31 & 1.28 & 0.107 & 1.31 & 1.60 & 0.110 & 1.31 & 1.42 \\
			            $32\times32$ & 0.238 & 1.31 & 1.30 & 0.186 & 1.31 & 1.66 & 0.111 & 1.31 & 1.47 \\
                        $64\times64$ & 0.240 & 1.33 & 1.32 & 0.249 & 1.33 & 1.73 & 0.115 & 1.33 & 1.52 \\
			          $128\times128$ & 0.239 & 1.33 & 1.32 & 0.279 & 1.33 & 1.75 & 0.114 & 1.33 & 1.53 \\
			\bottomrule
    \end{tabular}\label{tab:aniso_FE_newAMGr_improved_iter_scaling_crelax_cb3-gamgr}
\end{table}

\subsection{Algebraic Coarsening}\label{ssec:algebraic}

Given the satisfactory results in~\cref{tab:aniso_FE_newAMGr_improved_iter_scaling_crelax_cb3-gamgr}, we
next focus on extending these results to fully algebraic coarsening
using simulated annealing coarsening.  We note that this coarsening is computationally quite expensive~\cite{zaman2021coarse}, but that it provides the best known complexities for AMGr-style methods.  In~\cref{ssec:multilevel-gamgr}, we experiment with the more feasible greedy coarsening algorithm of MacLachlan and Saad~\cite{maclachlan2007greedy}. We emphasize that practical computing requires alternatives to the simulated annealing coarsening considered here, but that the greedy coarsening provides reasonable results at much more feasible cost.

{\bf Algebraic coarsening:} \Cref{tab:aniso_FE_newAMGr_improved_iter_scaling_crelax_sa_065-gamgr,tab:aniso_FE_newAMGr_improved_iter_scaling_crelax_sa_075-gamgr} show the convergence factors and corresponding grid and operator complexities for two-grid cycles using two values of the diagonal dominance parameter, $\eta$.  From preliminary experiments (not reported here), we noted substantial improvement when computing the fine-coarse partitioning using $\hat{A}$ (compared with $A$); hence, we use $\hat{A}$ in all subsequent results.

\begin{figure}
\centering
\includegraphics{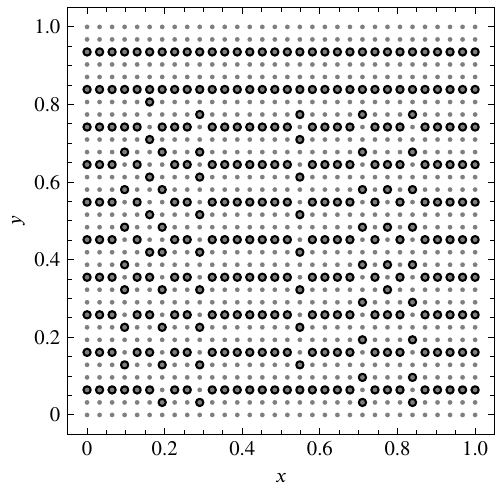}  %
\qquad
\includegraphics{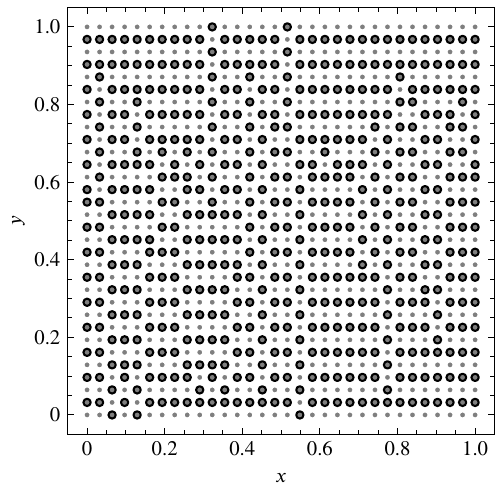}  %
\caption{Algebraic coarse-fine partitionings for the FE discretization of anisotropic-diffusion with $\theta = 0$ on uniform $32\times 32$ grid.  At left, partitioning with $\eta=0.65$.  At right, partitioning with $\eta=0.75$.  Fine-grid DoFs (points in $\Omega$) are denoted by filled grey dots; those that are in $C$ are marked with black circles.}\label{fig:C-F-nodes_aniso_struc_065_075-gamgr}
\end{figure}

\Cref{tab:aniso_FE_newAMGr_improved_iter_scaling_crelax_sa_065-gamgr} uses $\eta = 0.65$, resulting in two-level grid complexities, $C_{\rm{grid}}\approx 4/3$, matching that of the geometric semi-coarsening by three used above.  Using a larger parameter, $\eta = 0.75$, in~\cref{tab:aniso_FE_newAMGr_improved_iter_scaling_crelax_sa_075-gamgr}, results in $C_{\rm{grid}}\approx 3/2$, consistent with geometric semi-coarsening by a factor of two.  The coarsening is visualized in~\cref{fig:C-F-nodes_aniso_struc_065_075-gamgr} for both cases and $\theta = 0$, demonstrating that, while the coarsening is still algebraic, it retains much of the geometric character of semi-coarsening.  As expected, using a larger value of $\eta$ leads to an improvement in convergence factors (since the coarse-grid correction is over a larger space), but also higher complexities.
In particular, for the case of $\theta = 0$, we maintain complexities similar to those of geometric multigrid for these problems, with $C_{\rm{op}} \approx C_{\rm{grid}}$, but we also see the typical increase in AMG operator complexity faster than grid complexity for $\theta = \pi/6$ and $\pi/4$, indicating increased density in the coarse-grid operators.  While the effective convergence factors, defined as $\rho^{1/C_{\rm{op}}}$, are lower for $\eta = 0.75$ than $\eta = 0.65$, we emphasize that these are only two-grid complexities, and denser coarse-grid matrices lead to even higher three-grid complexities in the results to follow.  Thus, we focus on the choice of $\eta = 0.65$ in the results below.
\begin{table}
  \caption{Two-level AMGr convergence factors and corresponding complexities for anisotropic FE discretization {\bf using simulated annealing coarsening with $\eta = 0.65$}. $F$-nodes are
  interpolated using the SPAI approximation to $\hat{A}_{\ssff}^{-1}\hat{A}_{\ssfc}$, postprocessed to exactly interpolate a relaxed vector, and the SPAI approximations to $\hat{A}_{\ssff}^{-1}$ and $\hat{A}_{\sscc}^{-1}$ are used for $FCF$-relaxation, where $\hat{A}$ is the ``lumped'' matrix of strong connections computed from $A$.}
	\centering
		\begin{tabular}{c|N\Cgrid\CopTwo|N\Cgrid\CopTwo|N\Cgrid\CopTwo}  %
			\toprule
			      & \multicolumn{3}{c|}{$\theta=0$} & \multicolumn{3}{c|}{$\theta=\pi/6$} & \multicolumn{3}{c}{$\theta=\pi/4$}\\
     Grid size & $\rho$ & $C_{\text{grid}}$ & $C_{\text{op}}$ & $\rho$ & $C_{\text{grid}}$ & $C_{\text{op}}$ & $\rho$ & $C_{\text{grid}}$ & $C_{\text{op}}$ \EndTableHeader\\
     \midrule
                        $16\times16$ & 0.219 & 1.32 & 1.32 & 0.110 & 1.35 & 1.69 & 0.116 & 1.30 & 1.40 \\
			            $32\times32$ & 0.235 & 1.32 & 1.33 & 0.189 & 1.35 & 1.78 & 0.114 & 1.32 & 1.49 \\
                        $64\times64$ & 0.234 & 1.34 & 1.37 & 0.342 & 1.36 & 1.91 & 0.121 & 1.33 & 1.54 \\
			          $128\times128$ & 0.231 & 1.34 & 1.39 & 0.400 & 1.36 & 1.96 & 0.133 & 1.34 & 1.57 \\
			\bottomrule
		\end{tabular}\label{tab:aniso_FE_newAMGr_improved_iter_scaling_crelax_sa_065-gamgr}
\end{table}

\begin{table}
  \caption{Two-level AMGr convergence factors and corresponding complexities for anisotropic FE discretization {\bf using simulated annealing coarsening with $\eta = 0.75$}. $F$-nodes are
  interpolated using the SPAI approximation to $\hat{A}_{\ssff}^{-1}\hat{A}_{\ssfc}$, postprocessed to exactly interpolate a relaxed vector, and the SPAI approximations to $\hat{A}_{\ssff}^{-1}$ and $\hat{A}_{\sscc}^{-1}$ are used for $FCF$-relaxation, where $\hat{A}$ is the ``lumped'' matrix of strong connections computed from $A$.}
	\centering
		\begin{tabular}{c|N\Cgrid\CopTwo|N\Cgrid\CopTwo|N\Cgrid\CopTwo}  %
			\toprule
			      & \multicolumn{3}{c|}{$\theta=0$} & \multicolumn{3}{c|}{$\theta=\pi/6$} & \multicolumn{3}{c}{$\theta=\pi/4$}\\
     Grid size & $\rho$ & $C_{\text{grid}}$ & $C_{\text{op}}$ & $\rho$ & $C_{\text{grid}}$ & $C_{\text{op}}$ & $\rho$ & $C_{\text{grid}}$ & $C_{\text{op}}$ \EndTableHeader\\
     \midrule
                        $16\times16$ & 0.168 & 1.50 & 1.51 & 0.106 & 1.48 & 1.97 & 0.062 & 1.47 & 1.63 \\
			            $32\times32$ & 0.174 & 1.50 & 1.53 & 0.176 & 1.50 & 2.04 & 0.069 & 1.49 & 1.72 \\
                        $64\times64$ & 0.182 & 1.50 & 1.55 & 0.227 & 1.50 & 2.09 & 0.075 & 1.51 & 1.75 \\
			          $128\times128$ & 0.189 & 1.51 & 1.56 & 0.260 & 1.51 & 2.14 & 0.075 & 1.51 & 1.77 \\
			\bottomrule
		\end{tabular}\label{tab:aniso_FE_newAMGr_improved_iter_scaling_crelax_sa_075-gamgr}
\end{table}

{\bf Interpolation truncation:} To attenuate the higher complexities observed in~\cref{tab:aniso_FE_newAMGr_improved_iter_scaling_crelax_sa_065-gamgr,tab:aniso_FE_newAMGr_improved_iter_scaling_crelax_sa_075-gamgr}, we employ interpolation truncation, which is used successfully in several AMG settings~\cite{KStuben_2001a, HDeSterck_UMYang_JJHeys_2006a, HDeSterck_etal_2008b}.  In our approach, we first use SPAI (with a sparsity pattern of $\hat{A}_{\ssfc} + \hat{A}_{\ssff}\hat{A}_{\ssfc}$) to compute $W\approx \hat{A}_{\ssff}^{-1}\hat{A}_{\ssfc}$.  Then, we sweep row-wise through $W$, dropping entries that are less than a factor of $\zeta$ of the largest entry (by absolute value) in the row, to yield $\widehat{W}$.  As a final step, we compute matrix $S$ to match interpolation to the relaxed vector, $\vec{z}$, so that $\vec{z}_{\ssf} = S\widehat{W}\vec{z}_{\ssc}$.  \Cref{tab:aniso_FE_newAMGr_improved_iter_scaling_crelax_trunc_1_5th_sa_065-gamgr} shows results using $\zeta = 0.2$.  For the cases of $\theta = 0$ and $\pi/4$, we see that this truncation has no real effect in comparison with results in~\cref{tab:aniso_FE_newAMGr_improved_iter_scaling_crelax_sa_065-gamgr}.  This is easily understood from the nature of strong connections in these matrices, with only two strong connections per row, so there are only two interpolation weights in a typical row, and these weights are roughly equal in size.  In such cases, no effects of this truncation are expected.  For $\theta = \pi/6$, we see that this truncation leads to slight improvements in operator complexities, with small effects on convergence factors.  While the savings here may be minimal, we show below that interpolation truncation is an important tool in other cases, such as the isotropic Poisson problem in~\cref{sec:iso-poisson-gamgr}.  We also note that further increasing the truncation parameter to $\zeta = 0.25$ starts to show significantly degraded performance for $\theta = \pi/6$, when ``too many'' connections in interpolation are truncated.
\begin{table}
  \caption{Two-level AMGr convergence factors and corresponding complexities for anisotropic FE discretization using simulated annealing coarsening with $\eta = 0.65$. $F$-nodes are
  interpolated using the SPAI approximation to $\hat{A}_{\ssff}^{-1}\hat{A}_{\ssfc}$, postprocessed to exactly interpolate a relaxed vector, and the SPAI approximations to $\hat{A}_{\ssff}^{-1}$ and $\hat{A}_{\sscc}^{-1}$ are used for $FCF$-relaxation, where $\hat{A}$ is the ``lumped'' matrix of strong connections computed from $A$. {\bf Interpolation truncation with $\zeta = 0.2$ is used}.}
	\centering
		\begin{tabular}{c|N\Cgrid\CopTwo|N\Cgrid\CopTwo|N\Cgrid\CopTwo}  %
			\toprule
			      & \multicolumn{3}{c|}{$\theta=0$} & \multicolumn{3}{c|}{$\theta=\pi/6$} & \multicolumn{3}{c}{$\theta=\pi/4$}\\
     Grid size & $\rho$ & $C_{\text{grid}}$ & $C_{\text{op}}$ & $\rho$ & $C_{\text{grid}}$ & $C_{\text{op}}$ & $\rho$ & $C_{\text{grid}}$ & $C_{\text{op}}$ \EndTableHeader\\
     \midrule
                        $16\times16$ & 0.219 & 1.32 & 1.32 & 0.107 & 1.35 & 1.68 & 0.116 & 1.30 & 1.40 \\
			            $32\times32$ & 0.235 & 1.32 & 1.33 & 0.194 & 1.35 & 1.77 & 0.114 & 1.32 & 1.49 \\
                        $64\times64$ & 0.234 & 1.34 & 1.37 & 0.364 & 1.36 & 1.88 & 0.121 & 1.33 & 1.54 \\
			          $128\times128$ & 0.231 & 1.34 & 1.39 & 0.408 & 1.36 & 1.93 & 0.133 & 1.34 & 1.57 \\
			\bottomrule
		\end{tabular}\label{tab:aniso_FE_newAMGr_improved_iter_scaling_crelax_trunc_1_5th_sa_065-gamgr}
\end{table}

\Cref{tab:aniso_FE_newAMGr_improved_iter_crelax_trunc_1_5th_sa_065_3lev_V_Wcyc-gamgr} presents results using three-level cycles with $\zeta=0.2$ for these problems. As expected, adding more levels to the hierarchy increases the grid and operator complexities.  Using two-level geometric semi-coarsening-by-threes as a reference, we expect to see $C_{\rm{grid}} \approx 1 + \frac{1}{3} + \frac{1}{9} \approx 1.44$, which we do in all cases (with a slight increase for $\theta=\pi/6$, undoubtedly due to increased density of the coarse-grid operators).  We note that convergence does degrade going from two-grid to three-grid cycles, particularly for V-cycles (with convergence factors denoted by $\rho_V$), but also for W-cycles (with convergence factors denoted by $\rho_W$).%

\begin{table}
  \caption{{\bf Three-level} AMGr convergence factors and corresponding complexities for anisotropic FE discretization using simulated annealing coarsening with $\eta = 0.65$. $F$-nodes are
  interpolated using the SPAI approximation to $\hat{A}_{\ssff}^{-1}\hat{A}_{\ssfc}$, postprocessed to exactly interpolate a relaxed vector, and the SPAI approximations to $\hat{A}_{\ssff}^{-1}$ and $\hat{A}_{\sscc}^{-1}$ are used for $FCF$-relaxation, where $\hat{A}$ is the ``lumped'' matrix of strong connections computed from $A$. {\bf Interpolation truncation with $\zeta = 0.2$ is employed}.}
	\centering
		\begin{tabular}{c|NN\Cgrid\CopThree|NN\Cgrid\CopThree|NN\Cgrid\CopThree}  %
			\toprule
			     & \multicolumn{4}{c|}{$\theta=0$} & \multicolumn{4}{c|}{$\theta=\pi/6$} & \multicolumn{4}{c}{$\theta=\pi/4$}\\
                 Grid size & $\rho_V$ & $\rho_W$ & $C_{\text{grid}}$ & $C_{\text{op}}$ & $\rho_V$ & $\rho_W$ & $C_{\text{grid}}$ & $C_{\text{op}}$ & $\rho_V$ & $\rho_W$ & $C_{\text{grid}}$ & $C_{\text{op}}$ \EndTableHeader\\
     \midrule
                       $16\times16$ & 0.238 & 0.220 & 1.40 & 1.36 & 0.175 & 0.113 & 1.46 & 1.85 & 0.153 & 0.116 & 1.41 & 1.57 \\
			           $32\times32$ & 0.248 & 0.235 & 1.42 & 1.41 & 0.359 & 0.231 & 1.47 & 2.06 & 0.323 & 0.114 & 1.44 & 1.70 \\
                       $64\times64$ & 0.243 & 0.234 & 1.45 & 1.49 & 0.556 & 0.411 & 1.49 & 2.24 & 0.437 & 0.215 & 1.46 & 1.79 \\
			         $128\times128$ & 0.262 & 0.231 & 1.46 & 1.52 & 0.670 & 0.542 & 1.50 & 2.31 & 0.534 & 0.298 & 1.47 & 1.87  \\
			\bottomrule
		\end{tabular}\label{tab:aniso_FE_newAMGr_improved_iter_crelax_trunc_1_5th_sa_065_3lev_V_Wcyc-gamgr}
\end{table}

\subsection{The Generalized AMGr algorithm}\label{ssec:algorithm-gamgr}

Before presenting more extensive numerical results, we summarize the outcome of the experiments above in algorithmic form.  \Cref{alg:GenAMGrConst-gamgr} presents the AMGr setup algorithm: we compute the lumped matrix, $\hat{A}$, after finding strong connections, construct the $F$-$C$ partitioning based on $\hat{A}$, then form interpolation and the Galerkin coarse-grid operator.  In addition, we compute SPAI approximations to the inverses of $\hat{A}_{\ssff}$ and $\hat{A}_{\sscc}$, for use in relaxation, using $\mathcal{S}(A)$ to denote the sparsity pattern of matrix $A$.  As in the original AMGr paper, we use weighted relaxation, with optimal weights for single-step relaxation on the $F$ and $C$ subproblems, computed based on eigenvalue estimates for the ``preconditioned'' matrices, $M_{\ssff}A_{\ssff}$ and  $M_{\sscc}A_{\sscc}$.  For the coarse-grid problem, the computed eigenvalues are very close to one in all cases, and replacing the weighted relaxation with unweighted relaxation (approximating $\sigma_{\ssc} =1$) has little effect on convergence.  For the fine-grid relaxation, we find more variation in $\lambda_{\rm min}$ and $\lambda_{\rm max}$; we explore more practical alternatives to determining relaxation weights in~\cref{ssec:multilevel-gamgr}, although note that it may also be possible to choose better weights row-wise (e.g., by using weighting similar to the $\ell^1$-Jacobi relaxation method~\cite{doi:10.1137/100798806} or the matrix $S$ determined for interpolation).

\begin{algorithm}[t]
\caption{Generalized AMGr Setup Phase}\label{alg:GenAMGrConst-gamgr}
\begin{algorithmic}[1]
\Function{gen-amgr-setup}{$A, b$}

\State $\hat{A} \gets$ lumped approximation to $A$ after removing weak connections
\State $C, F \gets $ partitioning based on $\hat{A}$
\State $\hat{A}_{\ssff}, \hat{A}_{\ssfc}, \hat{A}_{\sscc} \gets$ extract submatrices of $\hat{A}$ based on $F$ and $C$

\State $P \gets \textsc{interpolation}(\hat{A}, \hat{A}_{\ssff}, \hat{A}_{\ssfc}, F, C)$
\State $A_{\ssc} \gets P^T A P$

\State $M_{\ssff} \gets \textsc{spai}(\hat{A}_{\ssff}, I_{\ssff}, \mathcal{S}\left(\hat{A}_{\ssff}\right))$
\State $\lambda_{\rm min}, \lambda_{\rm max} \gets$ minimum and maximum eigenvalues of $M_{\ssff}A_{\ssff}$\label{line:eigvals_DffinvAff-gamgr}
\State $\sigma_{\ssf} \gets 2/(\lambda_{\rm min} + \lambda_{\rm max})$\label{line:sigma_f_comput-gamgr}
\State $M_{\sscc} \gets \textsc{spai}(\hat{A}_{\sscc}, I_{\sscc}, \mathcal{S}\left(\hat{A}_{\sscc}\right))$
\State $\lambda_{\rm min}, \lambda_{\rm max} \gets$ minimum and maximum eigenvalues of $M_{\sscc}A_{\sscc}$\label{line:eigvals_DccinvAcc-gamgr}
\State $\sigma_{\ssc} \gets 2/(\lambda_{\rm min} + \lambda_{\rm max})$\label{line:sigma_c_comput-gamgr}
\State \textbf{return} $F, C, P, A_{\ssc}, M_{\ssff}, \sigma_{\ssf}, M_{\sscc}, \sigma_{\ssc}$

\EndFunction%
\end{algorithmic}
\end{algorithm}

The interpolation operator is computed following~\cref{alg:Interpolation-gamgr}.  First, the nonzero pattern is determined for interpolation, followed by a SPAI approximation to $\hat{A}_{\ssff}^{-1}\hat{A}_{\ssfc}$ for this pattern.  Small entries may be truncated in $W$ at this stage in order to reduce complexity of the resulting cycle.  After truncation, we rescale interpolation row-wise, either using constant scaling, as in~\cref{alg:constant-iteration-gamgr}, or the improved iteration scaling, as in~\cref{alg:improved-iteration-gamgr}.

\begin{algorithm}[t]
\caption{Interpolation Operator}\label{alg:Interpolation-gamgr}
\begin{algorithmic}[1]
\Function{interpolation}{$\hat{A}$, $\hat{A}_{\ssff}$, $\hat{A}_{\ssfc}$, $F$, $C$}
\State $\mathcal{Z} \gets \mathcal{S}\left(\hat{A}_{\ssfc} + \hat{A}_{\ssff}\hat{A}_{\ssfc}\right)$
\State $W \gets \textsc{spai}(\hat{A}_{\ssff},\hat{A}_{\ssfc},\mathcal{Z})$
\State $W \gets W$ after truncating small entries
\If{improved iteration scaling to be used}
\State $W \gets \textsc{improved-iteration-scaling}(\hat{A}, W, F, C)$
\Else%
\State $W \gets \textsc{constant-scaling}(W, F, C)$
\EndIf%
\State $P \gets \begin{bmatrix} W \\ I \end{bmatrix}$

\State \textbf{return} $P$

\EndFunction%
\end{algorithmic}
\end{algorithm}

\begin{algorithm}[t]
\caption{Constant Scaling}\label{alg:constant-iteration-gamgr}
\begin{algorithmic}[1]
\Function{constant-scaling}{$W, F, C$}
\State $\vec{s} \gets \vec{1}_{\ssf} / (W\vec{1}_{\ssc})$ \Comment{Componentwise division}
\State $S \gets$ matrix with $\vec{s}$ on diagonal and elsewhere $0$
\State $W \gets SW$
\State \textbf{return} $W$
\EndFunction%
\end{algorithmic}
\end{algorithm}

\begin{algorithm}[t]
\caption{Improved Iteration Scaling}\label{alg:improved-iteration-gamgr}
\begin{algorithmic}[1]
\Function{improved-iteration-scaling}{$\hat{A}, W, F, C$}
\State $n_{\rm wj} \gets 5$, $\omega \gets 2/3$, $D \gets$ diagonal of $\hat{A}$
\State $\vec{z} \gets \vec{1}$ %
\For {$i \gets 1,\ldots,n_{\rm wj}$}
    \State $\vec{z} \gets (I-\omega D^{-1}\hat{A})\vec{z}$
\EndFor%
\State $\vec{s} \gets \vec{z}_{\ssf} / (W\vec{z}_{\ssc})$ \Comment{Componentwise division}
\State $S \gets$ matrix with $\vec{s}$ on diagonal and elsewhere $0$
\State $W \gets SW$
\State \textbf{return} $W$
\EndFunction%
\end{algorithmic}
\end{algorithm}

Finally, we present the two-level AMGr solution phase in~\cref{alg:GenAMGrSol-gamgr}.  This includes either $F$- or $FCF$-relaxation (or more general relaxation) both before and after the coarse-grid correction phase, as well as a standard Galerkin coarse-grid correction. In what follows, we use $FCF$ relaxation consistently in all results. While we present the algorithm without implementation details, we note that the algorithms here can be implemented in either the ``natural'' ordering of matrix $A$, or in the ``permuted'' ordering given in~\cref{eq:permuted-gamgr}.  In many ways, it is simpler to implement the algorithm after permuting $A$ into its $F$-$C$ ordering.

\begin{algorithm}[t]
\caption{Generalized AMGr Solution Phase}\label{alg:GenAMGrSol-gamgr}
\begin{algorithmic}[1]
\Function{two-level}{$A, \vec{b}, \vec{u}, F, C, P, A_{\ssc}, M_{\ssff}, \sigma_{\ssf}, M_{\sscc}, \sigma_{\ssc}$}

        \State $\vec{u} \gets$ $F$-relaxation on $A\vec{u}=\vec{b}$
    \If{C-relaxation to be used}
        \For {$j \gets 1,\ldots,\nu_1$}
            \State $\vec{u} \gets$ $C$-relaxation on $A\vec{u}=\vec{b}$
            \State $\vec{u} \gets$ $F$-relaxation on $A\vec{u}=\vec{b}$
        \EndFor
    \EndIf%
\State $\vec{r}_{\ssc} \gets P^{T}(\vec{b}-A\vec{u})$
\State $\vec{e}_{\ssc} \gets$ solution of $A_{\ssc}\vec{e}_{\ssc}=\vec{r}_{\ssc}$ \Comment{use direct solve}
\State $\vec{u} \gets \vec{u}+P\vec{e}_{\ssc}$
    \State $\vec{u} \gets$ $F$-relaxation on $A\vec{u}=\vec{b}$
    \If{C-relaxation to be used}
        \For {$j \gets 1,\ldots,\nu_2$}
            \State $\vec{u} \gets$ $C$-relaxation on $A\vec{u}=\vec{b}$
            \State $\vec{u} \gets$ $F$-relaxation on $A\vec{u}=\vec{b}$
        \EndFor
    \EndIf%
\State \textbf{return} $\vec{u}$
\EndFunction%
\end{algorithmic}
\end{algorithm}

\section{Results}\label{sec:results-gamgr}

While the algorithms given above were derived by focusing on
performance for finite-element discretizations of anisotropic
diffusion equations on uniform grids, we emphasize in this section
that this methodology appears to have much wider applicability.  Here,
we first evaluate the approach on isotropic diffusion equations, on
both structured and unstructured grids.  Furthermore, we consider
results for a classic ``four-quadrant'' problem, with piecewise
constant diffusion and reaction coefficients on a uniform grid, and
for constant-coefficient anisotropic diffusion on an unstructured
grid.  In all results before~\cref{ssec:multilevel-gamgr}, we use the simulated annealing coarsening algorithm described above with $\eta = 0.65$.

\subsection{Isotropic Poisson problem}\label{sec:iso-poisson-gamgr}

In this section, we consider the isotropic diffusion equation $-\Delta u = f$, with Dirichlet boundary conditions, first on uniform meshes of the unit square domain.  As a benchmark, the first block column in~\cref{tab:iso_FE_newAMGr_improved_iter_crelax_sa_065-gamgr} presents convergence for the classical AMGr algorithm using a diagonal approximation, $M_{\ssff}$, to $A_{\ssff}$ in both relaxation and interpolation with only $F$-relaxation.  We note that using $\eta = 0.65$ already yields a positive effect on convergence; using $\eta = 0.56$ (as considered in past work) leads to convergence factors around 0.7, instead of 0.37.  In either case, while the convergence factors are bounded away from unity independently of grid size, the convergence is suboptimal for AMG on the model Poisson equation on a uniform grid.  The remaining columns of~\cref{tab:iso_FE_newAMGr_improved_iter_crelax_sa_065-gamgr} present results for the algorithm of~\cref{ssec:algorithm-gamgr}, demonstrating substantial improvement in two-grid convergence and reasonable three-level convergence.  We also note that the two-level generalized AMGr algorithm using $F$-relaxation in place of $FCF$-relaxation also offers reasonable convergence factors of about 0.1.  Here, we see that while the grid complexities for these cycles are relatively reasonable, the operator complexities are high, above 3.0 for most of the three-level cycles.

\begin{table}
  \caption{AMGr convergence factors and complexities for {\bf isotropic FE discretization on structured grids}.  Results for classical (diagonal $M_{\ssff}$) AMGr appear in the first block column, followed by those for the two-level and three-level SPAI-based algorithm in~\cref{ssec:algorithm-gamgr} in following block columns, using $\zeta = 0$.}
	\centering
		\begin{tabular}{c|N\Cgrid\CopTwo|N\Cgrid\CopTwo|NN\Cgrid\CopThree}   %
			\toprule
			      & \multicolumn{3}{c|}{Classical Two-level cycle} & \multicolumn{3}{c|}{Two-level cycle} & \multicolumn{4}{c}{Three-level cycles}\\
     Grid size & $\rho$ & $C_{\text{grid}}$ & $C_{\text{op}}$ & $\rho$ & $C_{\text{grid}}$ & $C_{\text{op}}$ & $\rho_V$ & $\rho_W$ & $C_{\text{grid}}$ & $C_{\text{op}}$ \EndTableHeader\\
     \midrule
                        $16\times16$ & 0.365 & 1.36 & 1.71 & 0.041 & 1.36 & 2.22 & 0.097 & 0.044 & 1.49 & 2.46 \\
			            $32\times32$ & 0.375 & 1.38 & 1.80 & 0.041 & 1.38 & 2.64 & 0.094 & 0.042 & 1.53 & 3.07 \\
                        $64\times64$ & 0.373 & 1.40 & 1.86 & 0.039 & 1.40 & 2.86 & 0.102 & 0.042 & 1.55 & 3.37 \\
			          $128\times128$ & 0.367 & 1.41 & 1.90 & 0.040 & 1.41 & 2.97 & 0.120 & 0.042 & 1.57 & 3.54 \\
			\bottomrule
		\end{tabular}\label{tab:iso_FE_newAMGr_improved_iter_crelax_sa_065-gamgr}
\end{table}

\begin{figure}
\centering
\includegraphics[width=\textwidth]{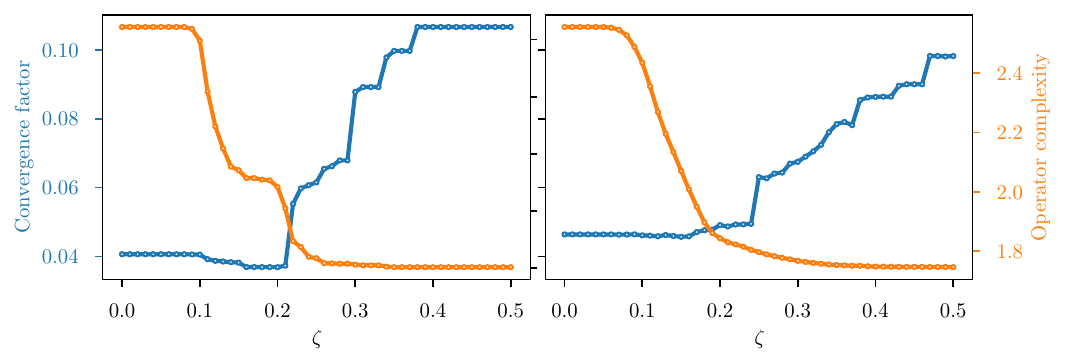}
\caption{Trade-off between two-level convergence factor and operator complexities as a function of $\zeta$, for isotropic Poisson on a uniform $32\times 32$ grid (at left) and on the unstructured triangulation with 1433 DoFs (at right).}\label{fig:convf_opcp_interp_trunc_iso_struc_unstruc_065-gamgr}
\end{figure}

To reduce the computational complexities, we truncate the smaller elements in the interpolation operator as discussed above.  A critical question in using interpolation truncation is the choice of the value of parameter $\zeta$.  The left plot of~\cref{fig:convf_opcp_interp_trunc_iso_struc_unstruc_065-gamgr} shows the effects of varying this parameter for the $32\times 32$ uniform grid.  For small values of $\zeta$, we observe large complexities, but also excellent two-level convergence factors.  As $\zeta$ increases past $0.2$, so do the convergence factors, yet the complexity continues drop.  If we were solely concerned with convergence, we might conclude that this is the optimal value of $\zeta$, since it yields the lowest complexity while retaining the best-possible convergence factor. However, to better balance cost vs.\ complexity, we prefer to take $\zeta =0.25$, where we approximately minimize the two-level complexity, while still retaining an acceptable convergence factor.  \Cref{tab:iso_FE_newAMGr_improved_iter_crelax_trunc_1_4th_sa_065-gamgr} shows two- and three-grid performance as we vary grid size with $\zeta = 0.25$. We see substantial improvements in complexity, with two-level complexities now similar to those of the classical AMGr algorithm in~\cref{tab:iso_FE_newAMGr_improved_iter_crelax_sa_065-gamgr}, and three-level grid complexities now about 2.2, instead of over 3.0.  At the same time, excellent two-level convergence factors are maintained, and there is only a slight impact on three-level convergence factors.

\begin{table}
  \caption{AMGr convergence factors and complexities for isotropic FE discretization on structured grids, using the two-level and three-level SPAI-based algorithm in~\cref{ssec:algorithm-gamgr}, {\bf with $\zeta = 0.25$}.}
	\centering
		\begin{tabular}{c|N\Cgrid\CopTwo|NN\Cgrid\CopThree}   %
			\toprule
			      & \multicolumn{3}{c|}{Two-level cycle} & \multicolumn{4}{c}{Three-level cycles}\\
     Grid size & $\rho$ & $C_{\text{grid}}$ & $C_{\text{op}}$ & $\rho_V$ & $\rho_W$ & $C_{\text{grid}}$ & $C_{\text{op}}$ \EndTableHeader\\
     \midrule
                        $16\times16$ & 0.053 & 1.36 & 1.73 & 0.099 & 0.058 & 1.49 & 1.90 \\
			            $32\times32$ & 0.061 & 1.38 & 1.84 & 0.140 & 0.071 & 1.52 & 2.07 \\
                        $64\times64$ & 0.065 & 1.40 & 1.91 & 0.151 & 0.074 & 1.55 & 2.20 \\
			          $128\times128$ & 0.069 & 1.41 & 1.95 & 0.157 & 0.079 & 1.56 & 2.26 \\
			\bottomrule
		\end{tabular}\label{tab:iso_FE_newAMGr_improved_iter_crelax_trunc_1_4th_sa_065-gamgr}
\end{table}

An important consideration for algebraic multigrid methods is whether or not they retain their performance as we transition from structured to unstructured grids.  Hence, our next problem considers the same isotropic diffusion operator, but discretized using piecewise linear finite elements on unstructured triangulations of the square domain, $[-1,1]^2$.  We construct grids by starting from an unstructured grid, performing several steps of uniform refinement, then smoothing the resulting grids.  Here, we consider three levels of refinement, generating meshes with \num{1433}, \num{5617}, and \num{22241} DoFs.  We again study the effects of varying the truncation parameter, $\zeta$, at right of~\Cref{fig:convf_opcp_interp_trunc_iso_struc_unstruc_065-gamgr}, and conclude that taking $\zeta = 0.25$ again gives a good trade-off between convergence and complexity.  \Cref{tab:iso_unstruc_FE_newAMGr_improved_iter_scaling_crelax_sa_065-gamgr} shows the resulting two- and three-grid convergence factors and operator complexities for the new AMGr algorithm applied to these problems.  While the convergence factors are somewhat larger than those for the uniform-grid discretization, they remain acceptable for AMG convergence for an isotropic diffusion operator.  Furthermore, we see the efficacy of interpolation truncation in reducing the operator complexity while maintaining acceptable convergence factors.

\begin{table}
  \caption{AMGr convergence factors and complexities for isotropic FE discretization on {\bf unstructured grids}, using the two-level and three-level SPAI-based algorithm in~\cref{ssec:algorithm-gamgr}. Results in the left-most block column show two-level results with no interpolation truncation ($\zeta = 0$), while the other block columns show results with $\zeta = 0.25$.}
	\centering
		\begin{tabular}{c|N\Cgrid\CopTwo|N\Cgrid\CopTwo|NN\Cgrid\CopThree}   %
			\toprule
			      & \multicolumn{3}{c|}{without truncation} & \multicolumn{7}{c}{with interpolation truncation}\\
			      & \multicolumn{3}{c|}{Two-level cycle} & \multicolumn{3}{c|}{Two-level cycle} & \multicolumn{4}{c}{Three-level cycles}\\
     \#DoF  & $\rho$ & $C_{\text{grid}}$ & $C_{\text{op}}$ & $\rho$ & $C_{\text{grid}}$ & $C_{\text{op}}$ & $\rho_V$ & $\rho_W$ & $C_{\text{grid}}$ & $C_{\text{op}}$\EndTableHeader\\
     \midrule
			            1433 & 0.046 & 1.36 & 2.56 & 0.063 & 1.36 & 1.80 & 0.155 & 0.064 & 1.50 & 2.17 \\
			            5617 & 0.135 & 1.36 & 2.53 & 0.132 & 1.36 & 1.75 & 0.167 & 0.127 &  1.50 & 2.17 \\
			           22241 & 0.167 & 1.37 & 2.49 & 0.167 & 1.37 & 1.78 & 0.322 & 0.186 &  1.52 & 2.22 \\
			\bottomrule
		\end{tabular}\label{tab:iso_unstruc_FE_newAMGr_improved_iter_scaling_crelax_sa_065-gamgr}
\end{table}

\subsection{Four-quadrant problems}\label{ssec:fourquad-gamgr}

\begin{figure}[!ht]
\centering
\includegraphics{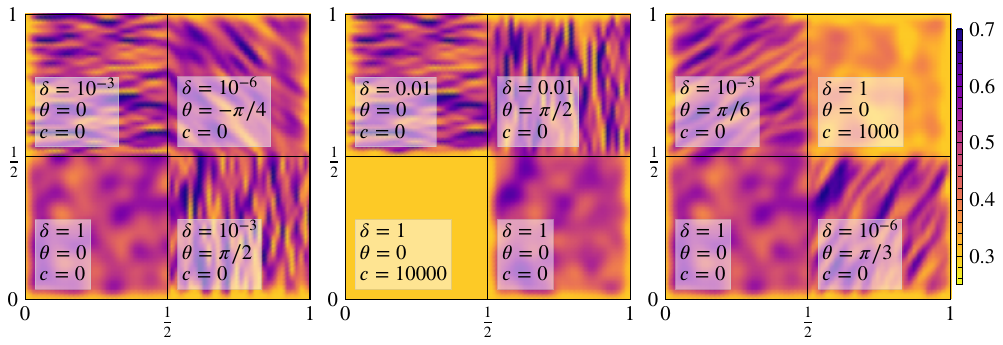}
\caption{Visualization of the strong connections for the four-quadrant problems; Problem 1 (left), Problem 2 (middle), and Problem 3 (right).}\label{fig:variable_coefficients_all_fourquad-gamgr}
\end{figure}

Next, we consider a family of two-dimensional anisotropic diffusion problems by adding a reaction term to Equation~\eqref{eq:anisotropic_diffusion-gamgr}, giving
\begin{equation}\label{eq:anisotropic_diffusion_reac-gamgr}
-\nabla \cdot \boldsymbol{K}(x, y)\nabla u(x, y) + c(x,y)u(x, y) = b(x,y)
\end{equation}
in the
domain $[0,1]\times [0,1]$ with Dirichlet boundary conditions.
The tensor coefficient is chosen as $\boldsymbol{K}(x, y)=QHQ^{T}$,
where $Q=\begin{bmatrix}
\cos(\theta) & -\sin(\theta)
\\ \sin(\theta) & \cos(\theta)
 \end{bmatrix}$, and $H=\begin{bmatrix}
1 & 0\\
0 & \delta
 \end{bmatrix}$, where $\theta$ specifies the
direction of anisotropy in the problem and $\delta$ specifies its strength. 
We partition the domain $[0,1]\times [0,1]$ into four equal quadrants and consider constant values of $\theta$, $\delta$ and $c$ within each quadrant, but with different values in different quadrants of the domain. The four-quadrant problem is common in AMG literature, and we consider three different problems within this class, with coefficient values shown below in \Cref{fig:variable_coefficients_all_fourquad-gamgr}.  Problem 1 is similar to the problem in Chapter 8 in the book by Briggs, Henson, and McCormick~\cite{briggs2000multigrid}, with no reaction term, large contrasts in the anisotropy strength, and non-grid-aligned diffusion in just one quadrant.  Problem 2 is the 2D-4Reg problem from Brannick and Falgout~\cite{JJBrannick_RDFalgout_2010a}, with a large reaction coefficient in one quadrant, but a small contrast in anisotropy strength and only grid-aligned anisotropy. Finally, Problem 3 is constructed to provide a more significant challenge, including a large reaction coefficient in one quadrant, a large contrast in anisotropy strength, and anisotropy directions in two quadrants that are neither aligned with the grid nor with the grid diagonal.

Two-level and three-level AMGr performance for these problems is shown in~\cref{tab:four_quad_newAMGr_improved_iter_crelax_trunc_1_4th_sa_065-gamgr,tab:four_quad_newAMGr_improved_iter_crelax_trunc_1_4th_sa_065_3lev_V_Wcyc-gamgr}, respectively.  We note that the two-level performance for Problems 1 and 2 is generally good, both in terms of convergence factor and complexity, while Problem 3 is clearly a harder problem.  Indeed, both V- and W-cycle convergence continues to perform well for Problem 2 in the three-level results in~\cref{tab:four_quad_newAMGr_improved_iter_crelax_trunc_1_4th_sa_065_3lev_V_Wcyc-gamgr}, with convergence outperforming that reported for the $65\times 65$ grid in Table 4.2 in Brannick and Falgout~\cite{JJBrannick_RDFalgout_2010a}, with comparable grid and operator complexities to the compatible relaxation AMG solver proposed there, and much better complexities than those reported there for BoomerAMG~\cite{VEHenson_UMYang_2002a}.  Problem 3 is clearly more taxing for AMG, yet the proposed generalized AMGr approach offers acceptable convergence in all cases.  Whether further improvement to these results is possible (or the performance degradation with grid size can be attenuated using Krylov acceleration) is left for future work.

\begin{table}
  \caption{{\bf Two-level} AMGr convergence factors and complexities for the {\bf four-quadrant problem} using the SPAI-based algorithm in~\cref{ssec:algorithm-gamgr}, with $\zeta = 0.25$.}
	\centering
		{\begin{tabular}{c|N\Cgrid\CopTwo|N\Cgrid\CopTwo|N\Cgrid\CopTwo}  %
			\toprule
			      & \multicolumn{3}{c|}{Problem 1} & \multicolumn{3}{c|}{Problem 2} & \multicolumn{3}{c}{Problem 3}\\
     Grid size & $\rho$ & $C_{\text{grid}}$ & $C_{\text{op}}$ & $\rho$ & $C_{\text{grid}}$ & $C_{\text{op}}$ & $\rho$ & $C_{\text{grid}}$ & $C_{\text{op}}$ \EndTableHeader\\
     \midrule
                        $17\times17$ & 0.388 & 1.35 & 1.49 & 0.602 & 1.26 & 1.32 & 0.062 & 1.25 & 1.38 \\
			            $33\times33$ & 0.433 & 1.34 & 1.51 & 0.581 & 1.26 & 1.37 & 0.266 & 1.33 & 1.62 \\
                        $65\times65$ & 0.413 & 1.36 & 1.56 & 0.595 & 1.27 & 1.41 & 0.491 & 1.39 & 1.78 \\
			          $129\times129$ & 0.420 & 1.36 & 1.57 & 0.599 & 1.36 & 1.60 & 0.692 & 1.39 & 1.80 \\

			\bottomrule
		\end{tabular}}\label{tab:four_quad_newAMGr_improved_iter_crelax_trunc_1_4th_sa_065-gamgr}
\end{table}

\begin{table}
  \caption{{\bf Three-level} AMGr convergence factors and complexities for the {\bf four-quadrant problem} using the SPAI-based algorithm in~\cref{ssec:algorithm-gamgr}, with $\zeta = 0.25$.}
	\centering
		{\begin{tabular}{c|NN\Cgrid\CopThree|NN\Cgrid\CopThree|NN\Cgrid\CopThree}  %
			\toprule
			        & \multicolumn{4}{c|}{Problem 1} & \multicolumn{4}{c|}{Problem 2} & \multicolumn{4}{c}{Problem 3}\\
           Grid size & $\rho_V$ & $\rho_W$ & $C_{\text{grid}}$ & $C_{\text{op}}$ & $\rho_V$ & $\rho_W$ & $C_{\text{grid}}$ & $C_{\text{op}}$ & $\rho_V$ & $\rho_W$ & $C_{\text{grid}}$ & $C_{\text{op}}$\EndTableHeader\\
     \midrule
                       $17\times17$ & 0.404 & 0.390 & 1.48 & 1.65 & 0.608 & 0.602 & 1.34 & 1.39 & 0.095 & 0.064 & 1.32 & 1.44 \\
			           $33\times33$ & 0.474 & 0.440 & 1.46 & 1.68 & 0.584 & 0.581 & 1.35 & 1.48 & 0.314 & 0.273 & 1.43 & 1.78 \\
                       $65\times65$ & 0.523 & 0.422 & 1.49 & 1.77 & 0.598 & 0.595 & 1.36 & 1.54 & 0.539 & 0.499 & 1.52 & 2.05 \\
			           $129\times129$ & 0.599 & 0.433 & 1.49 & 1.80 & 0.604 & 0.600 & 1.46 & 1.75 & 0.740 & 0.702 & 1.54 & 2.12 \\
			\bottomrule
		\end{tabular}}\label{tab:four_quad_newAMGr_improved_iter_crelax_trunc_1_4th_sa_065_3lev_V_Wcyc-gamgr}
\end{table}

\subsection{Unstructured anisotropic diffusion}\label{ssec:unstructured-gamgr}

\begin{figure}
\centering
\includegraphics[width=0.5\textwidth]{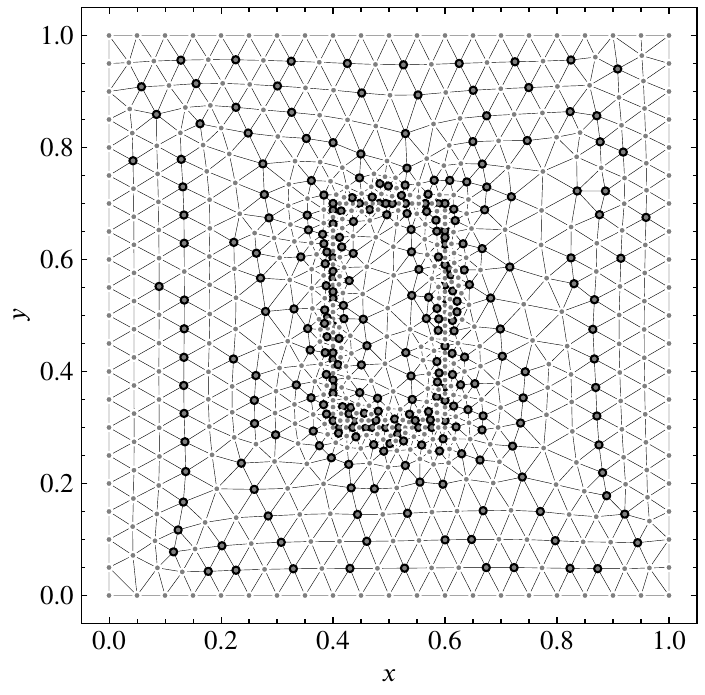}
\caption{Unstructured triangulation containing 798 points from Brannick and Falgout~\cite{JJBrannick_RDFalgout_2010a}, and its partitioning using $\eta = 0.65$. Fine-grid DoFs (points in $\Omega$) are denoted by filled grey dots; those that are in $C$ (\num{282} points) are marked with black circles.}\label{fig:cfnodes_unstruc_dom2_0p01_pi6_798_FE_lloyd_36_22_500000_1_ba_filt_lum_map_065dd-gamgr}
\end{figure}

Next, we consider the anisotropic diffusion problem in~\cref{eq:anisotropic_diffusion_reac-gamgr} with Dirichlet boundary conditions, $c=0$, $\theta=\pi/3$, and $\delta = 0.01$, on an unstructured triangulation of the unit square taken from Brannick and Falgout~\cite{JJBrannick_RDFalgout_2010a}, where the problem is labeled as 2D-M2-RLap. As in Table 4.2 from Brannick and Falgout~\cite{JJBrannick_RDFalgout_2010a}, we consider three refinements of the unstructured mesh for this problem, yielding discretized problems with \num{798}, \num{3109}, and \num{12273} DoFs, respectively. The mesh containing \num{798} DoFs (and its coarsening using $\eta = 0.65$) is shown in~\cref{fig:cfnodes_unstruc_dom2_0p01_pi6_798_FE_lloyd_36_22_500000_1_ba_filt_lum_map_065dd-gamgr}.  \Cref{tab:unstruc_dom2_0p01_pi6_FE_newAMGr_improved_iter_scaling_crelax_trunc_1_4th_sa_065-gamgr} presents two- and three-level convergence results for the generalized AMGr algorithm applied to this problem.  For comparison, we note that Table 4.2 of Brannick and Falgout~\cite{JJBrannick_RDFalgout_2010a} reports higher convergence factors (up to 0.95 on the finest grid) for CR-AMG applied to this problem, but at lower grid and operator complexities.  Compared to the BoomerAMG results presented in the same table, we see comparable convergence (0.57 for the finest grid) at lower complexity (2.6 at the finest grid, albeit for a multilevel cycle, not a three-level cycle).  Compared to classical AMGr applied to this problem, as given in Table 6 of Zaman et al.~\cite{zaman2021coarse}, we see substantial improvement in convergence factors (compared to values of 0.8--0.9 on the finest grid) and lower complexities in these results.

\begin{table}
  \caption{Two- and three-level AMGr convergence factors and complexities for the {\bf anisotropic diffusion problem on unstructured meshes} using the SPAI-based algorithm in~\cref{ssec:algorithm-gamgr}, with $\zeta = 0.25$.}
	\centering
		{\begin{tabular}{c|N\Cgrid\CopTwo|NN\Cgrid\CopThree}  %
			\toprule
			      & \multicolumn{3}{c|}{Two-level cycle} & \multicolumn{4}{c}{Three-level cycles} \\
     \#DoF  & $\rho$ & $C_{\text{grid}}$ & $C_{\text{op}}$ & $\rho_V$ & $\rho_W$ & $C_{\text{grid}}$ & $C_{\text{op}}$ \EndTableHeader\\
     \midrule
			             798 & 0.516 & 1.35 & 1.64 & 0.586 & 0.534 & 1.49 & 1.95 \\
			            3109 & 0.607 & 1.37 & 1.64 & 0.676 & 0.621 & 1.51 & 1.99 \\
			           12273 & 0.601 & 1.37 & 1.65 & 0.707 & 0.639 & 1.51 & 2.01 \\
			\bottomrule
		\end{tabular}}\label{tab:unstruc_dom2_0p01_pi6_FE_newAMGr_improved_iter_scaling_crelax_trunc_1_4th_sa_065-gamgr}
\end{table}

\subsection{Multilevel Results}\label{ssec:multilevel-gamgr}

Two major obstacles remain in~\cref{alg:GenAMGrConst-gamgr} for
transitioning from the two- and three-level cycles studied above to
standard multilevel cycles.  First of all, we have (until now) focused
on the use of simulated annealing for determining the partitioning of
$A$ and $\hat{A}$ into the $F$ and $C$ sets.  While this is very
effective, it is also very costly, as many SA steps are required to
generate near-optimal partitionings using this algorithm.  Thus, we
switch here to using the greedy coarsening algorithm of MacLachlan and
Saad~\cite{maclachlan2007greedy}, which is much more efficient, but
generates poorer-quality partitions.  To compensate, we
investigate the effect of the diagonal dominance parameter,
$\eta$, on the complexities and convergence of the resulting
multilevel hierarchies, in order to attenuate some of the complexity
growth that we observe in the initial results.

The second major obstacle is the calculation of extremal eigenvalues
in~\cref{line:eigvals_DffinvAff-gamgr,line:eigvals_DccinvAcc-gamgr}
of~\cref{alg:GenAMGrConst-gamgr}, which has additional heavy
computational cost.  To eliminate this, we replace the optimal
calculation of $\sigma_{\ssf}$ and $\sigma_{\ssc}$ with a common
heuristic estimate of the optimal regularization parameter.  Knowing
that $A_{\ssff}$ and $A_{\sscc}$ are both positive-definite matrices,
we expect that $M_{\ssff}$ and $M_{\sscc}$ are as
well.  If this is the case, the spectra of $M_{\ssff}A_{\ssff}$ and
$M_{\sscc}A_{\sscc}$ are guaranteed to be contained in the
intervals from 0 to their largest eigenvalues, which can be estimated
by their maximum absolute row sums (using Ger\v{s}gorin's theorem).
While the $F$-relaxation originally used in AMGr targets an optimal
reduction over all modes by estimating both ends of the spectrum of
$M_{\ssff}A_{\ssff}$, we propose a simpler heuristic of
choosing $\sigma_F$ and $\sigma_C$ to be $3/2$ divided by the maximum
absolute row sum of $M_{\ssff}A_{\ssff}$ and
$M_{\sscc}A_{\sscc}$, respectively.  The choice of weight $3/2$ in
this heuristic reflects the expectation that these matrices are
well-conditioned, so we need not use a weight as large as 2
(which would be optimal if we estimate the smallest eigenvalues as 0),
but that they are far from perfectly conditioned, so the weight should
be larger than 1.  Numerical tests confirm that using weight $3/2$ is
a good compromise~---~in some cases, some improvements are possible
with larger weights, but this leads to greatly degraded performance in
some cases as well.

In the results that follow, we emphasize measured convergence factors and complexities rather than timings.  Our solver codes are written in Python, using standard tools from the Numpy and Scipy libraries.  As such, certain routines are difficult to optimize in pure Python code.  We note, in particular, that the greedy coarsening code consistently requires over 50\% of the time for the setup phase, but that a compiled implementation using linked lists and an incomplete bin sort was used by MacLachlan and Saad~\cite{maclachlan2007greedy}, resulting in compute times similar to those for the classical AMG setup phase.  The SPAI implementation requires another 33\% of the total setup time, primarily due to inefficient sparse matrix and index arithmetic.  While we are unaware of an efficient Python implementation of any SPAI algorithm, recent literature (for example, the work of Anzt and co-authors~\cite{ANZT20181, 7836596}) reports performant implementations of incomplete sparse approximate inverse preconditioners that could be leveraged in its place.  Remaining significant setup costs again are due to inefficient matrix operations in Python, such as the lumping used to remove weak connections, that could also be ameliorated in compiled code.  Overall, this leads us to believe that a compiled implementation of the algorithm could be competitive with existing AMG implementations, although we leave such an implementation for future work.

As a comparison with the final three-level results in~\cref{tab:aniso_FE_newAMGr_improved_iter_crelax_trunc_1_5th_sa_065_3lev_V_Wcyc-gamgr},
\Cref{tab:aniso_FE_newAMGr_improved_iter_crelax_trunc_1_5th_sa_065_multilev_V_Wcyc_sigapprox-gamgr} shows convergence factors, complexities, and number of levels in the multigrid hierarchies ($n_{\rm l}$) for the multilevel algorithm.
These results show notable degradation in both operator and grid complexities, due to the use of greedy coarsening with $\eta = 0.65$ in contrast with the simulated annealing coarsening used in the previous results.  Nonetheless, we observe excellent W-cycle convergence factors in all cases (outperforming the earlier results for $\theta = 0$ and $\theta = \pi/6$), and consistent V-cycle convergence factors.
In experiments not reported here, we compared convergence to the case of using exact eigenvalues and found little difference in convergence overall.  Notably, when using the multigrid cycles as preconditioners for conjugate gradient, using the heuristic choice incurs at most 3 additional iterations over using cycles based on the exact eigenvalue computation.
\begin{table}
  \addtolength{\tabcolsep}{-2pt}
  \caption{{\bf Multi-level} AMGr convergence factors and
    corresponding complexities for anisotropic FE discretization using
    {\bf greedy coarsening} with $\eta = 0.65$. $F$-nodes are
  interpolated using the SPAI approximation to
  $\hat{A}_{\ssff}^{-1}\hat{A}_{\ssfc}$, postprocessed to exactly
  interpolate a relaxed vector, and the SPAI approximations to
  $\hat{A}_{\ssff}^{-1}$ and $\hat{A}_{\sscc}^{-1}$ are used for
  $FCF$-relaxation, where $\hat{A}$ is the ``lumped'' matrix of strong
  connections computed from $A$. Interpolation truncation with $\zeta
  = 0.2$ is employed.  {\bf Estimates of the eigenvalues are used in relaxation}.}
	\centering
		\begin{tabular}{c|N N\Cgrid\CopThree c|N N\Cgrid\CopThree c|N N\Cgrid\CopThree c}  %
			\toprule
			     & \multicolumn{5}{c|}{$\theta=0$} & \multicolumn{5}{c|}{$\theta=\pi/6$} & \multicolumn{5}{c}{$\theta=\pi/4$}\\
                 Grid size & $\rho_V$ & $\rho_W$ & $C_{\text{grid}}$ & $C_{\text{op}}$ & $n_{\rm l}$ & $\rho_V$ & $\rho_W$ & $C_{\text{grid}}$ & $C_{\text{op}}$ & $n_{\rm l}$ & $\rho_V$ & $\rho_W$ & $C_{\text{grid}}$ & $C_{\text{op}}$ & $n_{\rm l}$ \EndTableHeader\\
     \midrule
		$32\times32$   & 0.205 & 0.186 & 1.78 & 1.73 & 4 & 0.312 & 0.136 & 1.76 & 2.24 & 4 & 0.383 & 0.170 & 1.74 & 2.11 & 4\\
                $64\times64$   & 0.230 & 0.187 & 1.84 & 1.81 & 6 & 0.515 & 0.139 & 1.88 & 2.50 & 6 & 0.460 & 0.149 & 1.87 & 2.45 & 6\\
		$128\times128$ & 0.232 & 0.188 & 1.91 & 1.88 & 7 & 0.641 & 0.163 & 1.93 & 2.64 & 8 & 0.668 & 0.285 & 1.92 & 2.61 & 8\\
		$256\times256$ & 0.242 & 0.186 & 1.95 & 1.93 & 8 & 0.722 & 0.180 & 1.96 & 2.72 & 10 & 0.740 & 0.346 & 1.96 & 2.72 & 10\\
			\bottomrule
		\end{tabular}\label{tab:aniso_FE_newAMGr_improved_iter_crelax_trunc_1_5th_sa_065_multilev_V_Wcyc_sigapprox-gamgr}
\end{table}

Similar results are shown in~\cref{tab:fourquad_newAMGr_improved_iter_crelax_trunc_1_4th_sa_065_multilev_V_Wcyc_sigapprox-gamgr} for the four-quadrant problems from~\cref{ssec:fourquad-gamgr}, for comparison with~\cref{tab:four_quad_newAMGr_improved_iter_crelax_trunc_1_4th_sa_065_3lev_V_Wcyc-gamgr}.  Again, we note that the complexities are much higher than those reported earlier using the simulated annealing coarsening algorithm, but that this added complexity pays off in improved multilevel convergence. For Problem 2, we again compare to the results presented in Table 4.2 by Brannick and Falgout~\cite{JJBrannick_RDFalgout_2010a}, and see that this coarsening achieves comparable complexities to those reported there for $33\times 33$ and $65\times 65$ grids, but much better convergence.  Overall, we again see grid-independent W-cycle convergence for each problem, but growth in V-cycle convergence factors.  When run as preconditioners for CG, we find that W-cycles lead to convergence in 5--10 iterations (more for Problem 1, fewer for Problems 2 and 3), and V-cycle convergence in up to 14 iterations (again, with Problem 1 requiring most, and Problem 3 requiring fewest).
\begin{table}
  \addtolength{\tabcolsep}{-2pt}
  \caption{{\bf Multi-level} AMGr convergence factors and complexities for the {\bf four-quadrant problem} using the SPAI-based algorithm in~\cref{ssec:algorithm-gamgr}, with $\zeta = 0.25$.  {\bf Greedy coarsening and the heuristic eigenvalue estimates are used}.}
	\centering
		\begin{tabular}{c|N N\Cgrid\CopThree c|N N\Cgrid\CopThree c|N N\Cgrid\CopThree c}  %
			\toprule
			     & \multicolumn{5}{c|}{Problem 1} & \multicolumn{5}{c|}{Problem 2} & \multicolumn{5}{c}{Problem 3}\\
                 Grid size & $\rho_V$ & $\rho_W$ & $C_{\text{grid}}$ & $C_{\text{op}}$ & $n_{\rm l}$ & $\rho_V$ & $\rho_W$ & $C_{\text{grid}}$ & $C_{\text{op}}$ & $n_{\rm l}$ & $\rho_V$ & $\rho_W$ & $C_{\text{grid}}$ & $C_{\text{op}}$ & $n_{\rm l}$ \EndTableHeader\\
     \midrule
		$33\times33$   & 0.541 & 0.412 & 1.77 & 2.08 & 4 & 0.429 & 0.334 & 1.55 & 1.73 & 4 & 0.202 & 0.097 & 1.66 & 2.17 & 4\\
                $65\times65$   & 0.602 & 0.419 & 1.87 & 2.30 & 6 & 0.435 & 0.263 & 1.65 & 1.92 & 6 & 0.414 & 0.128 & 1.84 & 2.59 & 6\\
		$129\times129$ & 0.698 & 0.421 & 1.93 & 2.43 & 8 & 0.323 & 0.256 & 1.83 & 2.34 & 8 & 0.621 & 0.151 & 1.92 & 2.82 & 8\\
		$257\times257$ & 0.783 & 0.420 & 1.97 & 2.49 & 11 & 0.367 & 0.258 & 1.95 & 2.63 & 11 & 0.769 & 0.210 & 1.96 & 2.95 & 10\\
			\bottomrule
		\end{tabular}\label{tab:fourquad_newAMGr_improved_iter_crelax_trunc_1_4th_sa_065_multilev_V_Wcyc_sigapprox-gamgr}
\end{table}

As a final 2D test problem, we present results for the anisotropic diffusion problem on unstructured meshes considered in~\cref{ssec:unstructured-gamgr} and~\cref{tab:unstruc_dom2_0p01_pi6_FE_newAMGr_improved_iter_scaling_crelax_trunc_1_4th_sa_065-gamgr}.  Here, to explore the connection between the diagonal dominance parameter, $\eta$, and the resulting complexities and convergence factors, we consider $\eta = 0.65$ as before, along with $\eta = 0.60$ and $\eta = 0.56$.  \Cref{tab:unstruc_dom2_0p01_pi6_FE_newAMGr_improved_iter_scaling_crelax_trunc_1_4th_sa_065_multilev_V_Wcyc_sigapprox2-gamgr} shows that, as expected, complexities decrease and convergence factors generally increase as $\eta$ gets smaller, but that significant improvements in complexity are possible by using smaller $\eta$ without sacrificing substantial convergence.  In particular, comparing results for $\eta = 0.56$, we observe modest increases in grid complexity in comparison with those in~\cref{tab:unstruc_dom2_0p01_pi6_FE_newAMGr_improved_iter_scaling_crelax_trunc_1_4th_sa_065-gamgr}, possibly attributed to the increase from three-level to multi-level cycles.  Comparing these results to that presented in Table 4.2 by Brannick and Falgout~\cite{JJBrannick_RDFalgout_2010a}, we observe complexities better than those reported for BoomerAMG for these problems, albeit with slightly worse convergence, and slightly worse than those reported for compatible relaxation, but with better convergence.
\begin{table}
  \addtolength{\tabcolsep}{-2pt}
  \caption{{\bf Multi-level} AMGr convergence factors and complexities for {\bf anisotropic diffusion problem on unstructured meshes} using the SPAI-based algorithm in~\cref{ssec:algorithm-gamgr}, with $\zeta = 0.25$.  {\bf Greedy coarsening, for three values of $\eta$, and the heuristic eigenvalue estimates are used}.}
	\centering
		\begin{tabular}{c|N N\Cgrid\CopThree c|N N\Cgrid\CopThree c|N N\Cgrid\CopThree c}  %
			\toprule
			     & \multicolumn{5}{c|}{$\eta=0.56$} & \multicolumn{5}{c|}{$\eta=0.60$} & \multicolumn{5}{c}{$\eta=0.65$}\\
                 \#DoF & $\rho_V$ & $\rho_W$ & $C_{\text{grid}}$ & $C_{\text{op}}$ & $n_{\rm l}$ & $\rho_V$ & $\rho_W$ & $C_{\text{grid}}$ & $C_{\text{op}}$ & $n_{\rm l}$ & $\rho_V$ & $\rho_W$ & $C_{\text{grid}}$ & $C_{\text{op}}$ & $n_{\rm l}$ \EndTableHeader\\
     \midrule
798   & 0.612 & 0.560 & 1.54 & 2.06 & 4 & 0.601 & 0.529 & 1.66 & 2.37 & 4 & 0.537 & 0.467 & 1.76 & 2.56 & 4\\
3109  & 0.732 & 0.665 & 1.66 & 2.36 & 5 & 0.712 & 0.591 & 1.77 & 2.66 & 6 & 0.686 & 0.577 & 1.90 & 2.98 & 6\\
12273 & 0.777 & 0.653 & 1.71 & 2.51 & 7 & 0.762 & 0.630 & 1.83 & 2.87 & 8 & 0.742 & 0.576 & 1.98 & 3.33 & 9\\
			\bottomrule
		\end{tabular}\label{tab:unstruc_dom2_0p01_pi6_FE_newAMGr_improved_iter_scaling_crelax_trunc_1_4th_sa_065_multilev_V_Wcyc_sigapprox2-gamgr}
\end{table}

Finally, we apply the AMGr algorithm to three-dimensional model problems.  Here, we make two changes to the parameters considered above.  First, we take $\eta = 0.51$ in the results that follow, reflecting the decrease in natural diagonal dominance when moving from two- to three-dimensional problems.  Secondly, while the multilevel results reported above for 2D problems coarsen until the coarsest-grid size is less than 100 nodes, for 3D problems, we use a limit of 500 nodes, as coarser levels than this were seen to lead to convergence difficulties.  Table~\ref{tab:3d_iso_FE_newAMGr_improved_iter_crelax_trunc_1_4th_sa_051_multilev_V_Wcyc_sigapprox-gamgr} shows results for the three-dimensional Laplacian on the unit cube, discretized using linear finite-elements on ``uniform'' tetrahedral meshes, constructed by taking uniform hexahedral grids of the given size and the cutting each hexahedron into six tetrahedra.  While we see some growth in complexity compared to the two-dimensional case, these seem reasonable without further tuning of the algorithm to account for the change to 3D problems, particularly given the consistent convergence factors.  Results for two anisotropic diffusion problems are shown in Table~\ref{tab:3d_aniso_FE_newAMGr_improved_iter_crelax_trunc_1_4th_sa_051_multilev_V_Wcyc_sigapprox-gamgr}.  Here, we consider the problem
\begin{equation}\label{eq:3d_anisotropic_diffusion-gamgr}
-\nabla \cdot \boldsymbol{K}\nabla u(x, y, z) = b(x, y, z)
\end{equation}
on the unit cube domain, with diffusion tensor
\[
\boldsymbol{K} = R^T\begin{bmatrix} \delta & 0 & 0 \\ 0 & \delta & 0 \\ 0 & 0 & 1 \end{bmatrix}R\text{ with }R = \begin{bmatrix}
\cos(\theta) & 0 & \sin(\theta) \\ 0 & 1 & 0 \\ -\sin(\theta) & 0 & \cos(\theta)
  \end{bmatrix}\begin{bmatrix}
  \cos(\theta) & -\sin(\theta) & 0 \\
  \sin(\theta) & \cos(\theta) & 0 \\
  0 & 0 & 1
  \end{bmatrix}.
\]
We set $\delta = 10^{-6}$ and consider $\theta = 0$ and $\pi/6$.  Again, we note some increase in complexity from the two-dimensional case, but overall reasonable complexities and convergence factors.

\begin{table}
  \addtolength{\tabcolsep}{-2pt}
  \caption{{\bf Multi-level} AMGr convergence factors and
    complexities for {\bf 3D isotropic diffusion} using the SPAI-based algorithm in~\cref{ssec:algorithm-gamgr}, with $\zeta = 0.25$.  {\bf Greedy coarsening with $\eta = 0.51$ and the heuristic eigenvalue estimates are used}.}
	\centering
		\begin{tabular}{c|N N\Cgrid\CopThree c}  %
			\toprule
                 Grid size & $\rho_V$ & $\rho_W$ & $C_{\text{grid}}$ & $C_{\text{op}}$ & $n_{\rm l}$ \EndTableHeader\\
     \midrule
         $16\times16\times16$   & 0.134 & 0.116 & 1.37 & 1.99 & 3 \\
		 $32\times32\times32$   & 0.494 & 0.077 & 1.48 & 2.53 & 4 \\
         $64\times64\times64$   & 0.727 & 0.133 & 1.53 & 2.73 & 6 \\
			\bottomrule
		\end{tabular}\label{tab:3d_iso_FE_newAMGr_improved_iter_crelax_trunc_1_4th_sa_051_multilev_V_Wcyc_sigapprox-gamgr}
\end{table}

\begin{table}
  \addtolength{\tabcolsep}{-2pt}
  \caption{{\bf Multi-level} AMGr convergence factors and
    complexities for {\bf 3D anisotropic diffusion} using the SPAI-based algorithm in~\cref{ssec:algorithm-gamgr}, with $\zeta = 0.25$.  {\bf Greedy coarsening with $\eta = 0.51$ and the heuristic eigenvalue estimates are used}.}
	\centering
		\begin{tabular}{c|N N\Cgrid\CopThree c|N N\Cgrid\CopThree c}  %
			\toprule
			     & \multicolumn{5}{c|}{$\theta=0$} & \multicolumn{5}{c|}{$\theta=\pi/6$} \\
                 Grid size & $\rho_V$ & $\rho_W$ & $C_{\text{grid}}$ & $C_{\text{op}}$ & $n_{\rm l}$ & $\rho_V$ & $\rho_W$ & $C_{\text{grid}}$ & $C_{\text{op}}$ & $n_{\rm l}$  \EndTableHeader\\
     \midrule
         $16\times16\times16$   & 0.289 & 0.218 & 1.65 & 1.99 & 4 & 0.498 & 0.444 & 1.45 & 2.06 & 3  \\
	$32\times32\times32$	    & 0.307 & 0.215 & 1.77 & 2.29 & 5 & 0.704 & 0.494 & 1.59 & 2.71 & 5 \\
        $64\times64\times64$   & 0.324 & 0.217 & 1.84 & 2.49 & 6 & 0.824 & 0.525 & 1.67 & 3.17 & 6 \\
			\bottomrule
		\end{tabular}\label{tab:3d_aniso_FE_newAMGr_improved_iter_crelax_trunc_1_4th_sa_051_multilev_V_Wcyc_sigapprox-gamgr}
\end{table}

\section{Conclusions and future work}\label{sec:conclusion-gamgr}

Reduction-based AMG methods have been proposed and studied in many settings over the past 15 years, building effective solvers that can be more closely related to AMG convergence theory than many other heuristic methods.  In this paper, we aim to improve the practical performance of AMGr approaches by targeting tools that can greatly improve performance for anisotropic diffusion equations.  Through extensive numerical results, we show that the combination of using SPAI~\cite{grote1997parallel, doi:10.1137/S1064827503423500} to approximate $A_{\ssff}^{-1}$ along with tools to control sparsity leads to effective solvers for both anisotropic and isotropic diffusion operators on structured and unstructured grids.  In our view, this work points to weaknesses in the existing theory for AMGr-type methods, where approximations to $A_{\ssff}^{-1}A_{\ssfc}$ have not been considered (to our knowledge) in any context, providing an opportunity for future theoretical work to complete this and related~\cite{doi:10.1137/19M1256117,BUI2021114111} algorithmic development.  This may also provide new insights into desirable properties of SPAI-like approximations in this context. A key step in this work is to see if the good convergence properties observed here for W-cycles can also be robustly extended to V-cycles, as needed for effective parallel solvers. Additionally, further experiments are needed to see how to adapt the methodology proposed here to an even broader set of challenging problems, including the indefinite Helmholtz equation and convection-dominated flows.

\section*{Acknowledgments}

The work of S.P.M. was partially supported by an NSERC Discovery Grant. This work does not have any conflicts of interest.

\bibliography{refs}

\end{document}